\documentclass[12pt]{article}
\usepackage{a4,amsmath,amsthm}
\usepackage{eucal}
\usepackage{pb-diagram}
\usepackage{amsfonts}
\newcommand{\ar}{\renewcommand{\arraystretch}{1}}
\newcommand{\C}{\Bbb{C}}
\newcommand{\Z}{\Bbb{Z}}

\newcommand{\K}{\Bbb{K}}
\newcommand{\R}{\Bbb{R}}
\newcommand{\BH}{\Bbb{H}}

\DeclareMathOperator{\End}{End}
\DeclareMathOperator{\spin}{{\bf Spin}}
\DeclareMathOperator{\pin}{{\bf Pin}}

\DeclareMathOperator{\Id}{Id}

\DeclareMathOperator{\sech}{{\rm sech}}

\newcommand{\im}{\mbox{\rm Im}\,}
\newcommand{\Om}{{\bf\Omega}} 
\newcommand{\cA}{\mathcal{A}}

\newcommand{\cE}{\mathcal{E}}

\newcommand{\M}{{\bf\sf M}}

\newcommand{\sS}{{\sf S}}
\newcommand{\bi}{{\bf i}}
\newcommand{\bj}{{\bf j}}
\newcommand{\bk}{{\bf k}}

\newcommand{\bJ}{{\bf J}}
\newcommand{\bK}{{\bf K}}

\newcommand{\bS}{{\bf S}}

\newcommand{\mw}{{\scriptstyle\bigwedge}}

\newcommand{\cl}{C\kern -0.2em \ell}

\newcommand{\p}{\prime}
\newcommand{\e}{\mbox{\bf e}}
\newcommand{\bom}{\boldsymbol{\omega}}

\newtheorem{theorem}{Theorem}

\newtheorem{definition}{Definition}
\newtheorem{lem}{Lemma}
\newtheorem{cor}{Corollary}
\begin{document}
\title{Spinor Representations of Surfaces in 4--Dimensional 
Pseudo--Riemannian Manifolds}
\author{V.V. Varlamov\\
Department of Mathematics\\
Siberia State University of Industry\\
Novokuznetsk 654007, Russia\\
e-mail: root@varlamov.kemerovo.su}
\date{}
\maketitle
\begin{abstract}
Spinor representations of surfaces immersed into 4--dimensional 
pseudo--riemannian manifolds are defined in terms of minimal left ideals
and tensor decompositions of Clifford algebras. The classification of
spinor fields and Dirac operators on the immersed surfaces is given.
The Dirac--Hestenes spinor field on surfaces immersed into Lorentzian
manifolds and on surfaces conformally immersed into Minkowski spacetime
is defined.
\end{abstract}
\begin{description}
\item[Mathematics Subject Classification (1991):] 15A66, 53A50, 53A10, 35Q55
\item[Keywords:] Clifford algebras, spinors, minimal left ideals, spinor
bundles, immersed surfaces, Dirac operators, Weierstrass representation.
\end{description}
\section{Introduction}
One of the most interesting aspects of the relationship between differential
geometry of surfaces and Lie groups is a theory of spinor representations
of surfaces. This interest can be explained by two reasons: firstly, the
spinor representations of surfaces (SRS) are natural consequence of a
fundamental relation between differential geometry and Lie groups, secondly,
spinor representations as a whole are of great importance in theoretical
physics. On the other hand, nowadays theory of surfaces in itself is
finding increasing use in many areas of modern physics, above all in
soliton theory \cite{Sym85,Bob94,Bob99} and string theory \cite{GM93}.
\begin{sloppypar}
However, at the present time SRS are exhaustively studied only for
3--dimensional case (in particular, for the case of conformal immersions
of surfaces into a 3--dimensional euclidean space). Whereas a high--grade
by content (both mathematical and physical point of view) case of immersions
of surfaces into 4--dimensional manifolds still remains poorly studied, i.e.
a problem of the finding and classification of spinor representations of
surfaces immersed into 4d pseudo--riemannian manifolds is opened.
The present work is devoted to solving this problem.
\end{sloppypar}
Historically, spinor structures in the differentiable manifolds introduced
by Haefliger in 1956 \cite{Hae56} (see also \cite{BH,Mil63}). 
Further, spinor structures over the Riemann
surfaces are defined by Atiyah \cite{At71} (see also \cite{John80}). These
works 
were
served as a basis for the following construction of spinor representation
of minimal surfaces immersed into a 3--dimensional space \cite{Sul89}.
The Sullivan's results was generalized by Abresh onto a case of surfaces
with constant mean curvature \cite{Abr89}. It is well--known that conformal
immersions of the minimal surfaces into a 3--dimensional euclidean space
are described by a classical Weierstrass representation \cite{Weier}.
Thus, there is a close relationship between the spinor structures over
the Riemann surfaces (SRS) and conformal immersions (see \cite{KS95,KS96}).
Moreover, the Weierstrass representation has a natural formulation in terms
of SRS \cite{KS95}. The following important step in this direction was
maked by Taimanov \cite{Tai97a} (see also \cite{Tai97b,Tai97c,Tai98}). In
the work \cite{Tai97a} SRS was constructed on the basis of a so--called
generalized Weierstrass representation (GWR) which describes conformal
immersions of generic (non--minimal) surfaces into $\R^3$, and also a
globalization of this spinor representation was proposed by means of
introduction of a spinor bundle. GWR playing a key role in these works,
firstly, appears in \cite{Eisen} and further was rediscovered by Kenmotsu
\cite{Ken79} and Konopelchenko \cite{Kon1}. In the work \cite{Kon1}
a system of two differential equations (a so--called
2--dimensional Dirac equation) which coincides with a linear problem of 
a modified Veselov--Novikov hierarchy (mVN--hierarchy) has been considered
along with GWR.
Thus, there is
a relationship between the theory of conformal immersions of surfaces into
$\R^3$ and soliton theory, since integrable deformations of surfaces
are defined by the mVN--hierarchy. Moreover, it allows to express deformations
of spinor fields (smooth sections of the spinor bundles) via 
mVN--deformations \cite{Tai97a}.
The other important point in \cite{Tai97a} is establishing a relation
with the works of Hoffman and Osserman \cite{HO80,HO83,HO85} on the
generalized Gauss map. The following important work in the theory of SRS
is a Friedrich's paper \cite{Fr98}. The main advantage of \cite{Fr98} is a
consideration of SRS in the framework of a theory of the Dirac operator
on the spin manifolds (see \cite{Bau81,Bar91,Fr97,Amm98}).

As noted above, spinor representations of surfaces in 4d manifolds are not
studied in practice, however, at the present time in the papers of
Konopelchenko and Landolfi \cite{KL98a,Kon2,KonLan2}  a generalized
Weierstrass representation  has been considered
for surfaces conformally immersed into 4d
pseudo--euclidean spaces. The basic subject of these works is an extension
of the results obtained in \cite{Kon1,KT95,KT96} onto 4--dimensional spaces.
At this point GWR constructed on the basis of the generalized Gauss map
\cite{HO85}. In connection with this it should be noted that a Dirac
operator and SRS for surfaces conformally immersed into 4d complex space
are considered recently in the framework of the generalized Weierstrass
representation \cite{Var99a,Var99b}.

One of the main goals of the present research is a formulation of GWR for
surfaces immersed into 4d manifolds in terms of spinor bundles. An initial
point is an extension of Friedrich's results \cite{Fr98} onto 4d
pseudo--riemannian manifolds. According to widely accepted interpretation
a spinor field on the manifold is understood as a smooth section of the
spinor bundle. On the other hand, there exists a more profound definition
introduced by Chevalley \cite{Che54} and further developed in the works
\cite{Lou81,Cru87,Cru91}. It is a so--called {\it algebraic definition}
of the spinor field in which spinor is understood as an element of a minimal
left ideal of the Clifford algebra. The advantage of this definition is
obvious, since it allows to directly use basic facts and theorems of the
Clifford algebras theory at the study of the spin manifolds. 

The present paper is organized as follows. Since the Clifford algebras in
toto are the base of the spinor representations, then the basic facts
about these algebras are considered in the section 2. In the section 3
the algebraic definition of a spinor field on the surface immersed into
a 3--dimensional manifold is given in terms of a minimal left ideal of the
Pauli algebra. Further, a form of the spinor fields defined in the section 4
is a natural extension of the construction presented above (the section 3)
onto 4d pseudo--riemannian manifolds. At this point the Clifford algebras
of tangent bundles of 4d manifolds are quaternionic algebras, i.e. for any
4--dimensional Clifford algebra there exists a decomposition into the
tensor product of two quaternion algebras which further associated
respectively with tangent and normal bundles of the immersed surface.
The main goal of the section 5 is finding of a Dirac operator on the
immersed surface. In accordance with \cite{Bau89a,BFGK,Fr97} we suppose
that a spinor field on the ambient manifold is a real Killing spinor field.
The classification of the Dirac operators on the surfaces depends on a
metric of the ambient manifold and also on a metric of the immersed surface
(space--like and time--like surfaces). In virtue of a close relation with
theoretical physics an immersion of the surface into the Lorentzian
manifold is considered in more details. The relationship between spinor
fields on surfaces immersed into the Lorentzian manifold and a
Dirac--Hestenes spinor field \cite{Hest1,Hest67,Hest76} (which has an
important meaning in the electron theory \cite{Kel93,DR}) is established.
Further, in local consideration we have conformal immersions of surfaces
into 4d pseudo--euclidean spaces, which are considered in the section 6.\section{Algebraic Preliminaries}
In this section we will list some basic facts about Clifford algebras
which relevant for our further consideration.
Let $\K$ be a field of characteristic 0 $(\K=\R,\,\K=\Om,\,\K=\C)$, where
$\Om$ is a field of double numbers $(\Om=\R\oplus\R)$, and $\R,\,\C$ are
the fields of real and complex numbers, respectively. A Clifford algebra
over a field $\K$ is an algebra with $2^n$ basis elements: $\e_0$
(unit of the algebra), $\e_1,\e_2,\ldots,\e_n$ and the products of the
one--index elements $\e_{i_1i_2\ldots i_k}=\e_{i_1}\e_{i_2}\ldots\e_{i_k}$.
Over the field $\K=\R$ the Clifford algebra denoted as $\cl_{p,q}$, where
the indices $p,q$ correspond to the indices of the quadratic form
\[
Q=x^2_1+\ldots+x^2_p-\ldots-x^2_{p+q}
\]
of a vector space $V$ associated with $\cl_{p,q}$. A multiplication law
of $\cl_{p,q}$ defined by the following rule:
\begin{equation}\label{e1}
\e^2_i=\sigma(q-i)\e_0,\quad\e_i\e_j=-\e_j\e_i,
\end{equation}
where 
\begin{equation}\label{e2}
\sigma(n)=\begin{cases}
-1 & \text{if $n\leq 0$},\\
+1 & \text{if $n>0$}.
\end{cases}
\end{equation}
\begin{theorem}[{\rm Chevalley \cite{Che55}}] \label{tChe}
Let $V$ and $V^{\p}$ be vector spaces endowed with quadratic forms $Q$ and
$Q^{\p}$ over the field $\K$. Then a Clifford algebra 
$\cl(V\oplus V^{\p},\,Q\oplus Q^{\p})$ is naturally isomorphic to
$\cl(V,Q)\hat{\otimes}\cl(V^{\p},Q^{\p})$.
\end{theorem}
An important role in the theory of Clifford algebras played the square of
the volume element $\omega=\e_{12\ldots n}$, $n=p+q$:
\begin{equation}\label{e3}
\omega^2=\begin{cases}
-1 & \text{if $p-q\equiv 1,2,5,6\pmod{8}$},\\
+1 & \text{if $p-q\equiv 0,3,4,7\pmod{8}$}.
\end{cases}
\end{equation}
If $p+q$ is even and $\omega^2=1$, then $\cl_{p,q}$ is called {\it
positive} and respectively {\it negative} if $\omega^2=-1$. Or, in
accordance with (\ref{e3}):
\[
\begin{array}{ccc}
\cl_{p,q}>0 &\text{if}& p-q\equiv 0,4\pmod{8},\\
\cl_{p,q}<0 &\text{if}& p-q\equiv 2,6\pmod{8}.
\end{array}
\]
\begin{theorem}[{\rm Karoubi \cite[prop. 3.16]{Karo}}] \label{tKar}
1) If $\cl(V,Q)>0$, and $\dim V$ is even, then
\[
\cl(V\oplus V^{\p},Q\oplus Q^{\p})\simeq\cl(V,Q)\otimes\cl(V^{\p},Q^{\p}).
\]
2) If $\cl(V,Q)<0$, and $\dim V$ is even, then
\[
\cl(V\oplus V^{\p}, Q\oplus Q^{\p})\simeq\cl(V,Q)\otimes\cl(V^{\p},-Q^{\p}). 
\]
\end{theorem}
Further, let $\C_n=\C\otimes\cl_{p,q}$ and $\Om_{p,q}=\Om\otimes\cl_{p,q}$
be the Clifford algebras over the fields $\K=\C$ and $\K=\Om$, respectively.
\begin{theorem}[{\rm Rozenfel'd \cite{Roz55}}]\label{t1}
If $n=p+q$ is odd, then
\begin{eqnarray}
\cl_{p,q}&\simeq&\C_{p+q-1}\quad\text{if $p-q\equiv 1,5\pmod{8}$},\nonumber\\
\cl_{p,q}&\simeq&\Om_{p-1,q}\nonumber\\
&\simeq&\Om_{p,q-1}\quad\text{if $p-q\equiv 3,7\pmod{8}$}.\nonumber
\end{eqnarray}
\end{theorem}
\noindent
{\bf Example}. Let us consider the algebra $\cl_{0,3}$. According to the
theorem \ref{t1} we have $\cl_{0,3}\simeq\Om_{0,2}$, where $\Om_{0,2}$
is an algebra of elliptic biquaternions (it is a first so--called
Grassmann's extensive algebra introduced by Clifford in 1878 \cite{cliff}).
Since $\Om=\R\oplus\R$ and $\Om_{p,q}=\Om\otimes\cl_{p,q}$, we have
$\cl_{0,3}\simeq\cl_{0,2}\oplus\cl_{0,2}\simeq\BH\oplus\BH$, where $\BH$ is
a quaternion algebra.\\[0.4cm]
Generalizing this example we obtain
\begin{eqnarray}
\cl_{p,q}&\simeq&\cl_{p-1,q}\oplus\cl_{p-1,q}\nonumber\\
&\simeq&\cl_{p,q-1}\oplus\cl_{p,q-1}\quad\text{if $p-q\equiv 3,7\pmod{8}$}.
\label{e5'}
\end{eqnarray}

Over the field $\K=\C$ there is the analogous result \cite{Rash}
\begin{theorem}\label{t2}
When $p+q\equiv 1,3,5,7\pmod{8}$ the Clifford algebra over the field
$\K=\C$ decomposes into a direct sum of two subalgebras:
\[
\C_{p+q}\simeq\C_{p+q-1}\oplus\C_{p+q-1}.
\]
\end{theorem}

A minimal left (respectively right) ideal of $\cl_{p,q}$
is a set of type $I_{p,q}=\cl_{p,q}
e_{pq}$ (resp. $e_{pq}\cl_{p,q}$), where $e_{pq}$ is a primitive
idempotent, i.e., $e^2_{pq}=e_{pq}$ and $e_{pq}$ cannot be represented as a
sum of two orthogonal idempotents, i.e., $e_{pq}\neq f_{pq}+g_{pq}$, where
$f_{pq}g_{pq}=g_{pq}f_{pq}=0,\;f^2_{pq}=f_{pq},\,g^2_{pq}=g_{pq}$.
\begin{theorem}[{\rm Lounesto \cite{Lou81}}] \label{tLou}
A minimal left ideal of $\cl_{p,q}$ 
is of the
type $I_{p,q}=\cl_{p,q}e_{pq}$, where $e_{pq}=\frac{1}{2}(1+\e_{\alpha_1})
\ldots\frac{1}{2}(1+\e_{\alpha_k})$ is a primitive idempotent of $\cl_{p,q}$
and $\e_{\alpha_1},\ldots,\e_{\alpha_k}$ are commuting elements of the
canonical basis of $\cl_{p,q}$ such that $(\e_{\alpha_i})^2=1,\,(i=1,2,\ldots,
k)$ that generate a group of order $2^k,\;k=q-r_{q-p}$ and $r_i$ are the
Radon-Hurwitz numbers, defined by the recurrence formula $r_{i+8}=r_i+4$
and
\begin{center}{\rm
\begin{tabular}{lcccccccc}
$i$ & 0 & 1 & 2 & 3 & 4 & 5 & 6 & 7 \\ \hline
$r_i$ & 0 & 1 & 2 & 2 & 3 & 3 & 3 & 3
\end{tabular}.}
\end{center}
\end{theorem} 
\section{The Algebraic Definition of a Spinor Field on a 3D Manifold}
Let us consider a 3-dimensional oriented Riemannian manifold $M^{3,0}$
with a fixed spin structure and also an oriented surface $S^{2,0}$
isometrically immersed into $M^{3,0}$, $S^{2,0}\hookrightarrow M^{3,0}$.
At this point the surface $S^{2,0}$ is understood as a Riemannian
submanifold with some spin structure. Moreover, since the normal bundle
of the surface (precisely speaking, hypersurface) $S^{2,0}$ is trivial,
then the spin structure of $M^{3,0}$ induces a spin structure on the
surface $S^{2,0}\hookrightarrow M^{3,0}$. It is well--known that a spinor
field is a smooth section of the spinor bundle. Let $\Phi$ be a spinor
field on the manifold $M^{3,0}$.
It is obvious that the field
$\Phi$ is a section of 2--dimensional spinor bundle $S=Q\times_{\spin(3,0)}
\Delta_{3,0}$, since $\spin(3,0)\simeq\spin(0,3)$ and $\cl^+_{0,3}\simeq\BH$,
then $\dim\Delta_{3,0}=2$ (here $\BH$ is a quaternion algebra, $\cl^+_{p,q}$
is a Clifford algebra of all even elements). On the other hand,
$\spin(3,0)\simeq Sp(1)\simeq SU(2)\simeq S^3$ \cite{Port} and $\cl_{3,0}
\simeq\C_2\simeq\M_2(\C)$. Further, it takes to find a restriction of the
spinor bundle $S=Q\times_{\spin(3,0)}\Delta_{3,0}$ of the manifold $M^{3,0}$
onto a spinor bundle $S_{M^{2,0}}=Q\times_{\spin(2,0)}\Delta_{2,0}$ of the
surface $M^{2,0}$ conformally immersed into $M^{3,0}$. Let $\phi$ is a
spinor field on the surface $M^{2,0}$. Obviously, this field is a section
of 2--dimensional spinor bundle $S=Q\times_{\spin(2,0)}\Delta_{2,0}$, since
$\spin(2,0)\simeq\spin(0,2)$ and $\cl^+_{0,2}\simeq\cl_{0,1}\simeq\C$, then
$\spin(2,0)\simeq U(1)\simeq S^1$. Moreover, the spinor bundle of the surface
$M^{2,0}$ splits into two subbundles,
\[
S=S^+\oplus S^-,
\]
where $S^{\pm}=Q\times_{\spin(2,0)}\Delta^{\pm}_{2,0}$. Respectively,
a smooth section $\phi\in\Gamma(S)$ of the bundle $S$ has a form
$\varphi=\varphi^++\varphi^-$, where \cite{Fr98}
\begin{equation}\label{spin}
\varphi^+=\frac{1}{2}(\varphi+i\xi\cdot\varphi),\quad\varphi^-=
\frac{1}{2}(\varphi-i\xi\cdot\varphi),
\end{equation}
here $\xi=\e_3=\e_1\e_2$, $\varphi^+,\varphi^-$ are so--called half--spinors
(Weyl spinors) of the surface $S^{2,0}\hookrightarrow M^{3,0}$.

We will call the definition of the spinors (\ref{spin}) given above as
{\it a geometrical definition}, where the spinor field is understood as a
smooth section of the spinor bundle (this definition is widely used in
\cite{LM89,BFGK}). On the other hand, there exists {\it an algebraic definition}
of the spinor field as a minimal left ideal of the Clifford algebra
$\cl_{p,q}$ (see \cite{Che54,Lou81}). The algebraic definition in comparison
with geometric definition possess a more rich structure since allows
to directly use the all existing apparatus of the Clifford algebra theory.

Let us consider in details an algebraic definition of the spinors 
(\ref{spin}) as elements of a minimal left ideal $I_{3,0}=\cl_{3,0}e_{30}$
of the Pauli algebra $\cl_{3,0}$ ($\cl_{3,0}$ is a Clifford algebra of
a tangent bundle of the manifold $M^{3,0}$). In accordance with the
theorem \ref{tLou} a primitive idempotent of $\cl_{3,0}\simeq\C_2$ has a
form $e_{30}=\frac{1}{2}(1+\e_0)\sim\frac{1}{2}(1+i\e_{12})$, since in this
case a number of commuting elements equals to $k=q-r_{q-p}=0-r_{-3}=0-
(r_5-4)=1$. Further, it is obvious that $I_{3,0}=\cl_{3,0}e_{30}\simeq
\C_2e_{30}\simeq\M_2(\C)e_{30}$. By virtue of the isomorphism $\cl_{3,0}
\simeq\C_2\simeq\C\otimes\cl_{2,0}$ a general element of $\cl_{3,0}$
may be represented in the form of a following complex antiquaternion
\begin{equation}\label{complex}
\cA=\cl^0_{0,1}\e_0+\cl^1_{0,1}\e_1+\cl^2_{0,2}\e_2+\cl^3_{0,1}\e_1\e_2,
\end{equation}
where $\e^2_1=\e^2_2=1$. Since $\varphi\in I_{3,0}$, then
\begin{equation}\label{id1}
\begin{array}{ccc}
\varphi^+&=&\epsilon^+_{20}I_{3,0},\\
\varphi^-&=&\epsilon^-_{20}I_{3,0},
\end{array}
\end{equation}
where
\[
\epsilon^+_{20}=\frac{1}{2}(1+i\e_{12}),\quad
\epsilon^-_{20}=\frac{1}{2}(1-i\e_{12})
\]
are mutually orthogonal idempotents of the anti--quaternion (\ref{complex}).
Or, coming to matrix representations
\[
\ar\e_1\longmapsto\begin{pmatrix}
0 & 1\\
1 & 0
\end{pmatrix},\quad\e_2\longmapsto\begin{pmatrix}
0 & i\\
-i & 0
\end{pmatrix}
\]
we obtain that $\ar\varphi=\begin{pmatrix}\varphi_1 & 0 \\ \varphi_2 & 0
\end{pmatrix}\in I_{3,0}\simeq\M_2(\C)e_{30}$ and
\begin{eqnarray}
\varphi^+&=&\frac{1}{2}(1+i\e_{12})\varphi=\ar\begin{pmatrix}
\varphi_1 & 0\\
0 & 0
\end{pmatrix},\nonumber\\
\varphi^-&=&\frac{1}{2}(1-i\e_{12})\varphi=\ar\begin{pmatrix}
0 & 0\\
\varphi_2 & 0
\end{pmatrix},\label{id2}
\end{eqnarray} 
where
\begin{eqnarray}
\varphi_1&=&a^0-ia^{12},\nonumber\\
\varphi_2&=&a^1-ia^2,\nonumber
\end{eqnarray}
$a^0,a^1,a^2,a^{12}\in\C$. Thus, the Weyl spinors $\varphi^+,\varphi^-$
of the surface $S^{2,0}\hookrightarrow M^{3,0}$ are expressed via the
elements of the minimal left ideal of the Pauli algebra $\cl_{3,0}$ by the
formulae (\ref{id1}) or
(\ref{id2}).
\section{Spinor Structures on the Immersed Surfaces}
Let $M^{p,q}$ $(p+q=4)$ be a four--dimensional, real, connected, paracompact
manifold and let 
$TM^{p,q}$ (respectively $T^\ast M^{p,q}$) be a tangent (resp.
cotangent) bundle of the manifold $M^{p,q}$.
\begin{definition}
A Lorentzian manifold is a pair $(M^{1,3},g)$, where $g$ is a Lorentz
metric with a signature $(+,-,-,-)$, i.e. for any $x\in M^{1,3}$ there exists
an isomorphism $T_x M^{1,3}\simeq T^\ast_x M^{1,3}\simeq\R^{1,3}$, where 
$\R^{1,3}$ is a Minkowski spacetime.\\
A Majorana manifold is a pair $(M^{3,1},g)$, where $g$ is a metric with
a signature $(+,+,+,-)$, i.e. for any $x\in M^{3,1}$,
$T_x M^{3,1}\simeq T^\ast_x M^{3,1}\simeq\R^{3,1}$.\\
A Kleinian manifold is a pair $(M^{2,2},g)$, where a signature of the
metric $g$ has a form 
$(+,+,-,-)$ and for $x\in M^{2,2}$ follows $T_x M^{2,2}\simeq
T^\ast_x M^{2,2}\simeq\R^{2,2}$.
\end{definition}

Further, it is obvious that the spin structure of the 4--dimensional manifold
induces a spin structure of the immersed surface $S^{r,s}\hookrightarrow
M^{p,q}$, where $r+s=2,\,p+q=4$. In order to clarify this question let us
consider previously a more general case. Let
$M$ be an $(n+m)$--dimensional riemannian manifold and let 
$F\hookrightarrow M$ be an $n$--dimensional immersed submanifold.
We suppose that the both manifolds are endowed with a some spin structure.
Let $N$ be a normal bundle of the manifold $F\hookrightarrow M$.
In accordance with \cite{Mil65} the sum of the spin structures on the
tangent bundle and on the normal bundle of $F$ coincides with the spin
structure on the tangent bundle of $M$ restricted to $F$.
Thus, for any point $x\in F$ we have $T_xM=T_xF\oplus N_x$. Further, let
$\cl(T_xM)=\cl(T_xM,Q)$ be a Clifford algebra of the tangent space of the
manifold $M$ at the point $x$, where $Q$ is a quadratic form of a vector
space $\R^{p,q}\simeq T_xP$, $p+q=n+m$. In accordance with the 
theorem \ref{tChe} follows that $\cl(V\oplus V^{\p},Q\oplus Q^{\p})\simeq
\cl(V,Q)\hat{\otimes}\cl(V^{\p},Q^{\p})$, where in our case $V\simeq T_xF$,
$V^{\p}\simeq N_x$, $Q$ and $Q^{\p}$ are quadratic forms of the spaces
$T_xF$ and $N_x$. Moreover, if $\dim V$ is even we have the theorem
\ref{tKar}: $\cl(T_xF\oplus N_x,Q\oplus Q^{\p})
\simeq\cl(T_xF,Q)\otimes\cl(N_x,Q^{\p})$ if $\cl(T_xF,Q)>0$ and $\cl(T_xF\oplus
N_x,Q\oplus Q^{\p})\simeq\cl(T_xF,Q)\otimes\cl(N_x,-Q^{\p})$ if
$\cl(T_xF,Q)<0$. In any case the sum $T_xF\oplus N_x$ induces a tensor
product of the corresponding Clifford algebras.

Let us return to the 4--dimensional manifolds. It is easy to see that in this
case the sum $T_xF\oplus N_x$ induces a tensor product of the quaternion
algebras, $\cl_{r,s}\otimes\cl_{k,t}$, where $r+s=k+t=2$. At this point
there exist three types of tensor factors: a quaternion algebra $\cl_{0,2}$,
an anti--quaternion algebra $\cl_{2,0}$ and a pseudo--quaternion algebra 
$\cl_{1,1}$. Since a square of the volume element
$\omega=\e_{12}$ of $\cl_{0,2}$ equals to $-1$, then the algebra
$\cl_{0,2}$ is negative, $\cl_{0,2}<0$. Analogously, $\cl_{2,0}<0$
and $\cl_{1,1}>0$. Thus, according to the theorem \ref{tKar} we have the
following decompositions for the 4--dimensional Clifford algebras:
\begin{eqnarray}
\cl_{4,0}&\simeq&\cl_{2,0}\otimes\cl_{0,2},\nonumber\\
\cl_{0,4}&\simeq&\cl_{0,2}\otimes\cl_{2,0},\nonumber\\
\cl_{1,3}&\simeq&\cl_{1,1}\otimes\cl_{0,2},\nonumber\\
&\simeq&\cl_{0,2}\otimes\cl_{1,1},\nonumber\\
\cl_{3,1}&\simeq&\cl_{1,1}\otimes\cl_{2,0},\nonumber\\
&\simeq&\cl_{2,0}\otimes\cl_{1,1},\nonumber\\
\cl_{2,2}&\simeq&\cl_{2,0}\otimes\cl_{2,0},\nonumber\\
&\simeq&\cl_{0,2}\otimes\cl_{0,2},\nonumber\\
&\simeq&\cl_{1,1}\otimes\cl_{1,1}.\nonumber
\end{eqnarray}

Further, let $Q\times_{\spin(r,s)}\Delta_{r,s}$ be a spinor bundle of the
surface $S^{r,s}$, where $\Delta_{r,s}$ is a representation of the group
$\spin(r,s)$. Since $r+s=2$, then the bundle $S$ splits into two subbundles
$S=S^+\oplus S^-$, where $S^\pm=Q\times_{\spin(r,s)}\Delta^\pm_{r,s}$,
$\Delta^\pm_{r,s}$ are the spaces of half--spinors. The smooth section 
$\psi=\psi^+\oplus\psi^-\in\Gamma(S^+)\oplus\Gamma(S^-)$ is called
{\it a spinor field on the surface $S^{r,s}$}. On the other hand, the
components of the spinor field $\psi$ may be represented by elements of
the minimal left ideal of $\cl_{r,s}$. Let us consider in details the
minimal left ideals of the algebras $\cl_{2,0}\simeq\R(2),\,\cl_{0,2}\simeq
\BH$ and $\cl_{1,1}\simeq\R(2)$, which are the Clifford algebras of the
tangent spaces for three different types of the surfaces $S^{2,0},\,S^{0,2}$
and $S^{1,1}$. In accordance with the theorem \ref{tLou} we have
$k_{|\cl_{2,0}}=0-r_{-2}=0-(r_6-4)=1,\,k_{|\cl_{0,2}}=2-r_2=0,\,
k_{|\cl_{1,1}}=1-r_0=1$. Therefore, for the primitive idempotents of
$\cl_{2,0},\,\cl_{0,2}$ and $\cl_{1,1}$ we obtain respectively 
$e_{20}=\frac{1}{2}(1+\e_1),\,e_{02}=1$ and $e_{11}=\frac{1}{2}(1+\e_{12})$.
The minimal left ideals have a form $I_{2,0}=\cl_{2,0}\frac{1}{2}(1+\e_1),\,
I_{0,2}\equiv\cl_{0,2},\,I_{1,1}=\cl_{1,1}\frac{1}{2}(1+\e_{12})$ and
respectively for division rings we obtain $\K=e_{20}\cl_{2,0}e_{20}=\{1\},\,
\K=e_{02}\cl_{0,2}e_{02}=\{1,\e_1,\e_2,\e_{12}\},\,\K=e_{11}\cl_{1,1}e_{11}=
\{1\}$. Thus, on the surface $S^{2,0}$ there exists a spinor field
$\psi=\psi^++\psi^-$, where
\begin{eqnarray}
\psi^+&=&\epsilon^+_{20}I_{2,0},\nonumber\\
\psi^-&=&\epsilon^-_{20}I_{2,0},\nonumber
\end{eqnarray}
here
\[
\epsilon^+_{20}=\frac{1}{2}(1+\e_0),\quad\epsilon^-_{20}=\frac{1}{2}
(1-\e_0)
\]
are mutually orthogonal idempotents of $\cl_{2,0}$. Analogously, on the
surfaces $S^{0,2}$ and $S^{1,1}$ there exist spinor fields
\[
\psi^\pm=\epsilon^\pm_{02}I_{0,2},\quad\psi^\pm=\epsilon^\pm_{11}I_{1,1},
\]
where respectively
\[
\epsilon^\pm_{02}=\frac{1}{2}(1\pm i\e_{12}),\quad\epsilon^\pm_{11}=
\frac{1}{2}(1\pm\e_1).
\]

Further, let $\Psi$ be a spinor field on the 4--dimensional 
pseudo--Riemannian manifold. The our main goal is finding of restrictions
$\psi=\Psi_{|S^{r,s}}$ on the immersed surfaces $S^{r,s}\hookrightarrow
M^{p,q}$. First of all, the spinor field of the 4d manifold is a smooth
section of the spinor bundle, $\Psi=\Psi^++\Psi^-\in\Gamma(S^+)\oplus
\Gamma(S^-)$, where $S^\pm=Q\times_{\spin(p,q)}\Delta^\pm_{p,q}$,
$\Psi^\pm=\ar\begin{pmatrix}\Psi^\pm_1\\ \Psi^\pm_2\end{pmatrix}$.
By analogy with the 2--dimensional case the field $\Psi$ on $M^{p,q}$
may be considered as an element of the minimal left ideal of the
corresponding algebra $\cl_{p,q}$. According to the theorem \ref{tLou}
primitive idempotents of the algebras $\cl_{4,0},\,\cl_{0,4},\,\cl_{2,2},\,
\cl_{1,3},\,\cl_{3,1}$ have respectively a form: $e_{40}=\frac{1}{2}
(1+\e_1),\,e_{04}=\frac{1}{2}(1+\e_{123}),\,e_{22}=\frac{1}{2}(1+\e_{13})
\frac{1}{2}(1+\e_{24}),\,e_{31}=\frac{1}{2}(1+\e_{12})\frac{1}{2}(1+\e_1),\,
e_{13}=\frac{1}{2}(1+\e_{12})$. It should be noted that the different
structure of the presented idempotents may be explained by using of the
Karoubi theorem \ref{tKar}. Indeed, for the Majorana algebra $\cl_{3,1}$
by theorem \ref{tKar} there exists a decomposition $\cl_{3,1}\simeq
\cl_{1,1}\otimes\cl_{2,0}$ (or $\cl_{2,0}\otimes\cl_{1,1}$), which induces
a product of the primitive idempotents $e_{11}e_{20}$ (or $e_{20}e_{11}$).
Here $e_{11}=\frac{1}{2}(1+\e_{12})$ and $e_{20}=\frac{1}{2}(1+\e_1)$,
therefore $e_{31}=e_{11}e_{20}\sim\frac{1}{2}(1+\e_{24})\frac{1}{2}(1+\e_1)$.
In contrast with this, the spacetime algebra $\cl_{1,3}$ has a
decomposition $\cl_{1,3}\simeq\cl_{1,1}\otimes\cl_{0,2}$ ($\cl_{0,2}\otimes
\cl_{1,1}$), whence $e_{13}=e_{11}e_{02}\sim\frac{1}{2}(1+\e_{14})$, since
$e_{02}=1$. The analogous situation takes place for the algebras
$\cl_{4,0}\simeq\cl_{2,0}\otimes\cl_{0,2}$ and $\cl_{0,4}\simeq\cl_{0,2}
\otimes\cl_{2,0}$. Thus, for the minimal left ideals of the 4--dimensional
Clifford algebras and their division rings we have
\begin{eqnarray}
I_{0,4}&=&\cl_{0,4}\frac{1}{2}(1+\e_{123}),
\quad\K=\{1,\e_1,\e_{13},\e_3\}\simeq\BH;
\nonumber\\
I_{1,3}&=&\cl_{1,3}\frac{1}{2}(1+\e_{14}),\quad\K=\{1,\e_2,\e_3,\e_{23}\}
\simeq\BH;\nonumber\\
I_{2,2}&=&\cl_{2,2}\frac{1}{2}(1+\e_{13})\frac{1}{2}(1+\e_{24}),\quad
\K=\{1\}\simeq\R;\nonumber\\
I_{3,1}&=&\cl_{3,1}\frac{1}{2}(1+\e_1)\frac{1}{2}(1+\e_{24}),\quad
\K=\{1\}\simeq\R;\nonumber\\
I_{4,0}&=&\cl_{4,0}\frac{1}{2}(1+\e_1),\quad\K=\{1,\e_{23},\e_{24},\e_{34}\}
\simeq\BH.\nonumber
\end{eqnarray}

Let us consider now a spinor field of the time--like surface $S^{1,1}$
immersed into the Lorentzian manifold $M^{1,3}$. In this case the
Clifford algebra of a tangent space at the point $x\in S^{1,1}$ of the
manifold $M^{1,3}$ restricted to $S^{1,1}$ has a form $\cl_{1,3}\simeq
\cl_{1,1}\otimes\cl_{0,2}$. Let consider in more details a structure of
the decomposition $\cl_{1,1}\otimes\cl_{0,2}$. First of all, a similar form
of decomposition tells that for the algebra $\cl_{1,3}$ there exists
a transition from the real coordinates to quaternion coordinates of the
form $a+b\zeta_1+c\zeta_2+d\zeta_1\zeta_2$, where $\zeta_1=\e_{123},\,
\zeta_2=\e_{124}$ and $\zeta^2_1=\zeta^2_2=(\zeta_1\zeta_2)^2=-1$,
$\e^2=1,\e^2_2=\e^2_3=\e^2_4=-1$. The units $\zeta_1,\,\zeta_2$ are form
a basis of the quaternion algebra, since $\zeta_1\sim\bi,\,\zeta_2\sim\bj,\,
\zeta_1\zeta_2\sim\bk$. Therefore, a general element of the spacetime
algebra $\cl_{1,3}$ may be written as follows
\begin{equation}\label{quat}
\cA=\cl^0_{1,1}+\cl^1_{1,1}\zeta_1+\cl^2_{1,1}\zeta_2+\cl^3_{1,1}\zeta_1\zeta_2,
\end{equation}
where the every coefficient $\cl^i_{1,1}$ is isomorphic to the 
pseudo--quaternion algebra $\cl_{1,1}$:
\begin{eqnarray}
\cl^0_{1,1}&=&a+a^1\e_1+a^2\e_2+a^{12}\e_{12},\nonumber\\
\cl^1_{1,1}&=&-a^{123}-a^{23}\e_1-a^{13}\e_2-a^3\e_{12},\nonumber\\
\cl^2_{1,1}&=&-a^{124}-a^{24}\e_1+a^{14}\e_2+a^4\e_{12},\nonumber\\
\cl^3_{1,1}&=&-a^{34}-a^{134}\e_1-a^{234}\e_2+a^{1234}\e_{12}.\nonumber
\end{eqnarray}
It is easy to verify that the units $\zeta_i$ commute with every basis
element of the algebra $\cl_{1,1}$.

Further, let $\gamma:\,\cl_{1,3}\rightarrow\End(I_{1,3})$ be a spinor
representation of the spacetime algebra defined by the standard matrix
representation
\[
\ar\gamma_0=\gamma^0=\begin{pmatrix}
1 & 0\\
0 &-1
\end{pmatrix},\quad\gamma_k=-\gamma^k=\begin{pmatrix}
0 & -\sigma_k\\
\sigma_k & 0
\end{pmatrix},
\]
where $\sigma_k$ are the Pauli matrices
\[\ar
\sigma_1=\begin{pmatrix}
0 & 1\\
1 & 0
\end{pmatrix},\quad\sigma_2=\begin{pmatrix}
0 & -i\\
i & 0
\end{pmatrix},\quad\sigma_3=\begin{pmatrix}
1 & 0\\
0 & -1
\end{pmatrix}.
\]
By virtue of an identity $\cl_{1,3}\frac{1}{2}(1+\gamma_0)=\cl^+_{1,3}
\frac{1}{2}(1+\gamma_0)$ \cite{FRO90a} the minimal left ideal of $\cl_{1,3}$
takes a form $I_{1,3}=\cl^+_{1,3}\frac{1}{2}(1+\gamma_0)\simeq\cl_{3,0}
\frac{1}{2}(1+\gamma_0)$, since $\cl_{3,0}\simeq\cl^+_{1,3}$. Let
$\phi\in\cl_{3,0}$ be a Dirac--Hestenes spinor field and let
$\Psi\in I_{1,3}=\cl_{1,3}\frac{1}{2}(1+\gamma_0)$ be a so--called
{\it mother spinor} \cite{Lou93}, then
\[\ar
\Psi=\phi\frac{1}{2}(1+\gamma_0)=\begin{pmatrix}
\phi_1 & -\phi^\ast_2 & 0 & 0\\
\phi_2 & \phi^\ast_1 & 0 & 0\\
\phi_3 & \phi^\ast_4 & 0 & 0\\
\phi_4 & -\phi^\ast_3 & 0 & 0
\end{pmatrix},\quad\phi=\begin{pmatrix}
\phi_1 & -\phi^\ast_2 & \phi_3 & \phi^\ast_4\\
\phi_2 & \phi^\ast_1 & \phi_4 & -\phi^\ast_3\\
\phi_3 & \phi^\ast_4 & \phi_1 & -\phi^\ast_2\\
\phi_4 & -\phi^\ast_3 & \phi_2 & \phi^\ast_1
\end{pmatrix}.
\]
Let
\[
\epsilon^+_{13}=\frac{1}{2}(1+i\zeta_1\zeta_2),\quad\epsilon^-_{13}=\frac{1}{2}
(1-i\zeta_1\zeta_2)
\]
be mutually orthogonal idempotents of the quaternion (\ref{quat}), then by
analogy with (\ref{id1}) we have for a spinor field $\Psi$ of the
Lorentzian manifold $M^{1,3}$ restricted to
$S^{1,1}$ the following components
\begin{equation}\label{restrict}
\begin{array}{ccc}
\Psi^+_{|S^{1,1}}&=&\epsilon^+_{13}I_{1,3},\\
\Psi^-_{|S^{1,1}}&=&\epsilon^-_{13}I_{1,3},
\end{array}
\end{equation}
where $I_{1,3}\simeq\cl_{3,0}\frac{1}{2}(1+\gamma_0)\simeq
\M_2(\C)\frac{1}{2}(1+i\sigma_{12})$, since
$\cl_{3,0}\simeq\C_2\simeq\M_2(\C)$.
Therefore, for the mother spinor $\Psi\in I_{1,3}$ we have
\[
\Psi=\phi\frac{1}{2}(1+i\sigma_{12})=\ar\begin{pmatrix}
\phi^\ast_1+\phi^\ast_3 & 0\\
\phi_4-\phi_2 & 0
\end{pmatrix},
\]
where $\phi\in\C_2\simeq\M_2(\C)$ is a Dirac--Hestenes spinor field with
a following matrix representation
\[
\ar\phi=\begin{pmatrix}
\phi^\ast_1+\phi^\ast_3 & \phi^\ast_4+\phi^\ast_2\\
\phi_4-\phi_2 & \phi_1-\phi_3
\end{pmatrix},
\]
here
\[
\phi_1=a^0-ia^{12},\quad
\phi_2=-a^{13}-ia^{23},\quad
\phi_3=a^3-ia^{123},\quad
\phi_4=a^1+ia^2.
\]
Let us define matrix representations of the quaternion units 
$\zeta_1$ and $\zeta_2$ as follows
\[
\ar\zeta_1\;\longmapsto\;\begin{pmatrix}
0 & -1\\
1 & 0
\end{pmatrix},\quad\zeta_2\;\longmapsto\;\begin{pmatrix}
0 & i\\
i & 0
\end{pmatrix},
\]
then a restricted mother spinor $\Psi_{|S^{1,1}}$ on the time--like surface
$S^{1,1}\hookrightarrow M^{1,3}$ takes a form
\begin{eqnarray}
\Psi^+_{|S^{1,1}}&=&\frac{1}{2}(1+i\zeta_1\zeta_2)\Psi=\ar\begin{pmatrix}
\phi^\ast_1+\phi^\ast_3 & 0\\
0 & 0
\end{pmatrix},\nonumber\\
\Psi^-_{|S^{1,1}}&=&\frac{1}{2}(1-i\zeta_1\zeta_2)\Psi=\ar\begin{pmatrix}
0 & 0\\ 
\phi_4-\phi_2 & 0
\end{pmatrix}.\label{restrict2}
\end{eqnarray}

Analogously, for the space--like surface $S^{0,2}$ immersed into $M^{1,3}$
in virtue of the decomposition
$\cl_{1,3}\simeq\cl_{0,2}\otimes\cl_{1,1}$ the general element of
$\cl_{1,3}$ can be represented by a following pseudo--quaternion
\begin{equation}\label{pseudo}
\cA=\cl^0_{0,2}\zeta_0+\cl^1_{0,2}\zeta_1+\cl^2_{0,2}\zeta_2
+\cl^3_{0,2}\zeta_1\zeta_2,
\end{equation}
where $\zeta_1=\e_{134},\,\zeta_2=\e_{234}$ are pseudo--quaternion units,
$\zeta^2_1=-1,\,\zeta^2_2=1,\,
(\zeta_1\zeta_2)^2=1$, and an every coefficient $\cl^i_{0,2}$ in 
(\ref{pseudo}) is the quaternion algebra:
\begin{eqnarray}
\cl^0_{0,2}&=&a^0+a^3\e_3+a^4\e_4+a^{34}\e_{34},\nonumber\\
\cl^1_{0,2}&=&-a^{134}-a^{14}\e_3+a^{13}\e_4-a^1\e_{34},\nonumber\\
\cl^2_{0,2}&=&a^{234}+a^{24}\e_3-a^{23}\e_4-a^2\e_{34},\nonumber\\
\cl^3_{0,2}&=&a^{12}+a^{123}\e_3+a^{124}\e_4+a^{1234}\e_{34}.\nonumber
\end{eqnarray}
Let
\[
\ar\zeta_1\;\longmapsto\;\begin{pmatrix}
0 & -1\\
1 & 0
\end{pmatrix},\quad\zeta_2\;\longmapsto\;\begin{pmatrix}
0 & 1\\
1 & 0
\end{pmatrix},
\]
then for a restricted mother spinor $\Psi_{|S^{0,2}}$ on the space--like
surface $S^{0,2}\hookrightarrow M^{1,3}$ we obtain
\begin{eqnarray}
\Psi^+_{|S^{0,2}}&=&\frac{1}{2}(1+\zeta_1\zeta_2)I_{1,3},\nonumber\\
\Psi^-_{|S^{0,2}}&=&\frac{1}{2}(1-\zeta_1\zeta_2)I_{1,3},\nonumber
\end{eqnarray}
where $I_{1,3}\simeq\M_2(\C)\frac{1}{2}(1+i\sigma_{12})$.

Further, for the immersion $S^{1,1}\hookrightarrow M^{3,1}$ in virtue
of the decomposition $\cl_{3,1}\simeq\cl_{1,1}\otimes\cl_{2,0}$ a general
element of the Majorana algebra $\cl_{3,1}$ may be represented by a
following anti--quaternion
\begin{equation}\label{anti}
\cA=\cl^0_{1,1}+\cl^1_{1,1}\zeta_1+\cl^2_{1,1}\zeta_2+\cl^3_{1,1}\zeta_1
\zeta_2,
\end{equation}
where $\zeta_1=\e_{134},\,\zeta_2=\e_{234}$ are anti--quaternion units,
$\zeta^2_1=\zeta^2_2=1,\,(\zeta_1\zeta_2)^2=-1$, $\e^2_1=\e^2_2=\e^2_3=1,
\e^2_4=-1$. At this point an every coefficient $\cl^i_{1,1}$ in (\ref{anti})
is isomorphic to the pseudo--quaternion algebra:
\begin{eqnarray}
\cl^0_{1,1}&=&a^0+a^3\e_3+a^4\e_4+a^{34}\e_{34},\nonumber\\
\cl^1_{1,1}&=&a^{134}-a^{14}\e_3-a^{13}\e_4+a^1\e_{34},\nonumber\\
\cl^2_{1,1}&=&a^{234}-a^{24}\e_3-a^{23}\e_4+a^2\e_{34},\nonumber\\
\cl^3_{1,1}&=&-a^{12}+a^{123}\e_3+a^{124}\e_4+a^{1234}\e_{34}.\nonumber
\end{eqnarray}
Since the division ring of $\cl_{3,1}$ is $\K\simeq\R$, then $\cl_{3,1}
\simeq\R(4)$ and for a spinor representation of $\cl_{3,1}$ we have
$\gamma:\,\cl_{3,1}\rightarrow\End_{\R}(I_{3,1})$, where $I_{3,1}=\cl_{3,1}
e_{31}=\cl_{3,1}\frac{1}{2}(1+\e_1)\frac{1}{2}(1+\e_{24})$. The basis of
a spinspace $\sS\sim I_{3,1}$ defined as follows
\begin{eqnarray}
f_1&=&e_{31}=\frac{1}{4}(1+\e_1+\e_{24}+\e_{124}),\nonumber\\
f_2&=&\e_2e_{31}=\frac{1}{4}(\e_2-\e_{12}+\e_4-\e_{14}),\nonumber\\
f_3&=&\e_3e_{31}=\frac{1}{4}(\e_3-\e_{13}-\e_{234}+\e_{1234}),\nonumber\\
f_4&=&\e_{23}e_{31}=\frac{1}{4}(\e_{23}+\e_{123}-\e_{34}-\e_{134}).\nonumber
\end{eqnarray}
In this basis the matrices $\cE_i=\gamma(\e_i)$ are
\[\ar
\cE_1=\begin{pmatrix}
1 & 0 & 0 & 0\\
0 & -1 & 0 & 0\\
0 & 0 & -1 & 0\\
0 & 0 & 0 & 1
\end{pmatrix},\quad\cE_2=\begin{pmatrix}
0 & 1 & 0 & 0\\
1 & 0 & 0 & 0\\
0 & 0 & 0 & 1\\
0 & 0 & 1 & 0
\end{pmatrix},
\]
\begin{equation}\label{bas1}\ar
\cE_3=\begin{pmatrix}
0 & 0 & 1 & 0\\
0 & 0 & 0 &-1\\
1 & 0 & 0 & 0\\
0 &-1 & 0 & 0
\end{pmatrix},\quad\cE_4=\begin{pmatrix}
0 & 1 & 0 & 0\\
-1& 0 & 0 & 0\\
0 & 0 & 0 & 1\\
0 & 0 &-1 & 0
\end{pmatrix}.
\end{equation}
Let us show the validity of a relation $\cl_{3,1}\frac{1}{2}(1+\cE_1)
\frac{1}{2}(1+\cE_{24})\simeq\cl^+_{3,1}\frac{1}{2}(1+\cE_1)\frac{1}{2}
(1+\cE_{24})$ which plays a key role at the restriction of the bundle
$Q\times_{\spin(3,1)}\Delta_{3,1}$ of $M^{3,1}$ to the bundle
$Q\times_{\spin(1,1)}\Delta_{1,1}$ of the surface $S^{1,1}\hookrightarrow
M^{3,1}$. Since $\cl_{3,1}=\cl^+_{3,1}\oplus\cl^-_{3,1}$ and
$\cl_{3,1}\simeq\M_4(\R)$, then in the basis (\ref{bas1}) we obtain
\begin{equation}\label{rel}
\cl_{3,1}e_{31}=\cl^+_{3,1}e_{31}\oplus\cl^-_{3,1}e_{31}\simeq\ar
\begin{pmatrix}
0 & 0 & 0 & \phi_1\\
0 & 0 & 0 & \phi_2\\
0 & 0 & 0 & \phi_3\\
0 & 0 & 0 & \phi_4
\end{pmatrix}\oplus\begin{pmatrix}
0 & 0 & 0 & \eta_1\\
0 & 0 & 0 & \eta_2\\
0 & 0 & 0 & \eta_3\\
0 & 0 & 0 & \eta_4
\end{pmatrix},
\end{equation}
where  
\begin{gather}
\phi_1=-a^{23}+a^{34},\;\;\phi_2=a^{13}-a^{1234},\;\;\phi_3=-a^{12}-a^{14},\;\;
\phi_4=a^0+a^{24},\nonumber\\
\eta_1=-a^{123}+a^{134},\;\;\eta_2=-a^3+a^{234},\;\;\eta_3=a^2+a^4,\;\;
\eta_4=a^1+a^{124}.\nonumber
\end{gather}
If suppose $\xi_i=\phi_i+\eta_i$, then from (\ref{rel}) follows
$\cl_{3,1}e_{31}\simeq\cl^+_{3,1}e_{31}$. Further, let us define matrix
representations of the anti--quaternion units $\zeta_1$ and $\zeta_2$
as follows
\[\ar
\zeta_1\longmapsto\begin{pmatrix}
0 & 1\\
1 & 0
\end{pmatrix},\quad\zeta_2\longmapsto\begin{pmatrix}
0 & -i\\
i & 0
\end{pmatrix},
\]
then for components of the restricted spinor field of the time--like
surface $S^{1,1}\hookrightarrow M^{3,1}$ we obtain
\begin{eqnarray}
\Psi^+_{|S^{1,1}}&=&\frac{1}{2}(1+i\zeta_1\zeta_2)I_{3,1},\nonumber\\
\Psi^-_{|S^{1,1}}&=&\frac{1}{2}(1-i\zeta_1\zeta_2)I_{3,1},\label{immer}
\end{eqnarray}
where $I_{3,1}\simeq\cl^+_{3,1}\frac{1}{2}(1+\cE_1)\frac{1}{2}(1+\cE_{24})
\simeq\M_2(\C)\frac{1}{2}(1+i\sigma_{12})$ in virtue of an isomorphism
$\cl^+_{3,1}\simeq\cl_{1,2}\simeq\C_2$. Thus, as in the case of the
immersions $S^{1,1}\hookrightarrow M^{1,3},\,S^{0,2}\hookrightarrow M^{1,3}$
the components (\ref{immer}) expressed via the Dirac--Hestenes spinors
$\phi_i\in\C_2$.

The analogous restriction take place for an immersion $S^{2,0}\hookrightarrow
M^{3,1}$. In this case in virtue of the decomposition $\cl_{3,1}\simeq
\cl_{2,0}\otimes\cl_{1,1}$ the general element of $\cl_{3,1}$ represented
by a pseudo--quaternion $\cA=\sum^3_{i=0}\cl^i_{2,0}\zeta_i$, where
$\zeta_0=1,\zeta_1=\e_{124},\zeta_2=\e_{123},\zeta_3=\zeta_1\zeta_2$ and
$\zeta^2_1=1,\zeta^2_2=-1,\zeta^2_3=1$. At this point coefficients
$\cl^i_{2,0}$ (anti--quaternions) are
\begin{eqnarray}
\cl^0_{2,0}&=&a^0+a^1\e_1+a^2\e_2+a^{12}\e_{12},\nonumber\\
\cl^1_{2,0}&=&a^{124}+a^{24}\e_1-a^{14}\e_2-a^4\e_{12},\nonumber\\
\cl^2_{2,0}&=&a^{123}+a^{23}\e_1-a^{13}\e_2-a^3\e_{12},\nonumber\\
\cl^3_{2,0}&=&a^{34}+a^{134}\e_1+a^{234}\e_2+a^{1234}\e_{12}.\nonumber
\end{eqnarray}
The restricted spinor field on $S^{2,0}\hookrightarrow M^{3,1}$ has a
form $\Psi_{|S^{2,0}}=(\epsilon^+I_{3,1},\epsilon^-I_{3,1})$, where
$I_{3,1}\simeq\C_2\frac{1}{2}(1+i\sigma_{12}),\,\epsilon^\pm=\frac{1}{2}
(1\pm\zeta_3)$.

Let us consider now an immersion $S^{2,0}\hookrightarrow M^{4,0}$. First of
all, since the division ring of $\cl_{4,0}$ is $\K\simeq\BH$ we have
$\gamma:\,\cl_{4,0}\rightarrow\End_{\BH}(I_{4,0})$, where $I_{4,0}=\cl_{4,0}
\frac{1}{2}(1+\e_1)$, and for matrices $\cE_i=\gamma(\e_i)$ we obtain
respectively
\[\ar
\cE_1=\begin{pmatrix}
1 & 0 & 0 & 0\\
0 & 1 & 0 & 0\\
0 & 0 &-1 & 0\\
0 & 0 & 0 & -1
\end{pmatrix},\quad\cE_2=\begin{pmatrix}
0 & 0 & 1 & 0\\
0 & 0 & 0 & 1\\
1 & 0 & 0 & 0\\
0 & 1 & 0 & 0
\end{pmatrix},
\]
\begin{equation}\ar
\cE_3=\begin{pmatrix}
0 & 0 &-i & 0\\
0 & 0 & 0 &-i\\
i & 0 & 0 & 0\\
0 & i & 0 & 0
\end{pmatrix},\quad\cE_4=\begin{pmatrix}
0 & 0 & 0 & -i\\
0 & 0 & -i& 0\\
0 & i & 0 & 0\\
i & 0 & 0 & 0
\end{pmatrix}.\label{bas2}
\end{equation}
It is easy to verify that in the basis (\ref{bas2}) a relation
$\cl_{4,0}\frac{1}{2}(1+\cE_1)\simeq\cl^+_{4,0}\frac{1}{2}(1+\cE_1)$ holds.
Further, in virtue of $\cl_{4,0}\simeq\cl_{2,0}\otimes\cl_{0,2}$ a general
element of the algebra $\cl_{4,0}$ can be represented by a quaternion
$\cA=\sum^3_{i=0}\cl^i_{2,0}\zeta_i$, where $\zeta_1=\e_{123},\zeta_2=
\e_{124},\zeta_3=\zeta_1\zeta_2,\,\zeta^1_1=\zeta^2_2=\zeta^2_3=-1$.
Thus, for a restricted spinor field we have
\begin{eqnarray}
\Psi^+_{|S^{2,0}}&=&\frac{1}{2}(1+i\zeta_3)I_{4,0},\nonumber\\
\Psi^-_{|S^{2,0}}&=&\frac{1}{2}(1-i\zeta_3)I_{4,0},\nonumber
\end{eqnarray}
where $I_{4,0}\simeq\cl^+_{4,0}\frac{1}{2}(1+\cE_1)\simeq\Om_{0,2}\frac{1}{2}
(1+i\Upsilon_{12})$, since $\cl^+_{4,0}\simeq\cl_{0,3}\simeq\Om_{0,2}$
(theorem \ref{t1}), $\Om_{0,2}=\BH\oplus\BH$ is a semi--simple algebra
(an algebra of elliptic biquaternions), $\Upsilon_i$ are the matrix
representations of the units of $\Om_{0,2}$:
\[\ar
\Upsilon_1=\begin{pmatrix}
0 & -1\\
1 & 0
\end{pmatrix},\quad\Upsilon_2=\begin{pmatrix}
0 & i\\
i & 0
\end{pmatrix}.
\]
Let $\phi\in\Om_{0,2}$ be an 'elliptic' Dirac--Hestenes spinor field with
a following matrix representation
\[\ar
\phi=\begin{pmatrix}
\phi_1+\phi_3 & -\phi^\ast_4-\phi^\ast_2\\
\phi_4+\phi_2 & \phi^\ast_1+\phi^\ast_3
\end{pmatrix},
\]
where
\[
\phi_1=a^0-iea^{12},\quad\phi_2=ea^{23}+ia^2,\quad\phi_3=ea^{123}-ia^3,\quad
\phi_4=a^1+iea^{13},
\]
$e$ is a double unit. Then for components of the spinor field
$\Psi\in\Om_{0,2}\frac{1}{2}(1+i\Upsilon_{12})$ on the surface $S^{2,0}
\hookrightarrow M^{4,0}$ we obtain
\begin{eqnarray}
\Psi^+_{|S^{2,0}}&=&\frac{1}{2}(1+i\zeta_3)\Psi=\ar\begin{pmatrix}
\phi_1+\phi_3 & 0\\
0 & 0
\end{pmatrix},\nonumber\\
\Psi^-_{|S^{2,0}}&=&\frac{1}{2}(1-i\zeta_3)\Psi=\ar\begin{pmatrix}
0 & 0\\
\phi_4+\phi_2 & 0
\end{pmatrix}.\label{comp}
\end{eqnarray}
The immersion $S^{0,2}\hookrightarrow M^{0,4}$ is analogously defined.
In this case a general element of
$\cl_{0,4}$ is represented by an anti--quaternion $\sum^3_{i=0}\cl^i_{0,2}
\zeta_i$, where $\zeta_1=\e_{123},\zeta_2=\e_{124}$. Further, since
$\cl^+_{0,4}\simeq\cl_{0,3}\simeq\Om_{0,2}$, then $I_{0,4}\simeq\cl^+_{0,4}
\frac{1}{2}(1+\cE_{123})\simeq\Om_{0,2}\frac{1}{2}(1+i\Upsilon_{12})$.

Finally, let us consider spinor fields on the surfaces immersed into
the Kleinian manifold $M^{2,2}$. First of all, for a spinor representation
of the algebra $\cl_{2,2}$ we have $\gamma:\,\cl_{2,2}\rightarrow\End_{\R}
(I_{2,2})$, where $I_{2,2}=\cl_{2,2}\frac{1}{2}(1+\e_{13})\frac{1}{2}
(1+\e_{24})$, and for a basis of the spinspace $\sS\simeq I_{2,2}$ we have
also
\begin{eqnarray}
f_1&=&\frac{1}{4}(1+\e_{24}+\e_{13}-\e_{1234}),\nonumber\\
f_2&=&\frac{1}{4}(\e_1+\e_{124}+\e_3-\e_{234}),\nonumber\\
f_3&=&\frac{1}{4}(\e_2+\e_4-\e_{123}+\e_{134}),\nonumber\\
f_4&=&\frac{1}{4}(\e_{12}+\e_{14}-\e_{23}+\e_{34}).\nonumber
\end{eqnarray}
In this basis the matrices $\cE_i=\gamma(\e_i)$ are
\begin{gather}
\ar\cE_1=\begin{pmatrix}
0 & 1 & 0 & 0\\
1 & 0 & 0 & 0\\
0 & 0 & 0 &1\\
0 & 0 &1 & 0
\end{pmatrix},\quad\cE_2=\begin{pmatrix}
0 & 0 &1 & 0\\
0 & 0 & 0 &-1\\
1 & 0 & 0 & 0\\
0 & -1 & 0 & 0
\end{pmatrix},\nonumber\\
\cE_3=\ar\begin{pmatrix}
0 &-1 & 0 & 0\\
1 & 0 & 0 & 0\\
0 & 0 & 0 & -1\\
0 & 0 &1 & 0
\end{pmatrix},\quad\cE_4=\begin{pmatrix}
0 & 0 & -1 & 0\\
0 & 0 & 0 & 1\\
1 & 0 & 0 & 0\\
0 & -1 & 0 & 0
\end{pmatrix}.\label{bas3}
\end{gather}
It is easy to see that in the basis (\ref{bas3}) a relation
$\cl_{2,2}\frac{1}{2}(1+\cE_{13})\frac{1}{2}(1+\cE_{24})\simeq
\cl^+_{2,2}\frac{1}{2}(1+\cE_{13})\frac{1}{2}(1+\cE_{24})$ holds.
Further, for the immersion $S^{2,0}\hookrightarrow M^{2,2}$ a general
element of $\cl_{2,2}$ is represented by an anti--quaternion
$\sum^3_{i=0}\cl^i_{2,0}\zeta_i$, where $\zeta_1=\e_{123},\zeta_2=\e_{124}$,
and the coefficients $\cl^i_{2,0}$ (anti--quaternions) are generated by
a set $\{1,\e_1,\e_2,\e_{12}\}$.
By virtue of an isomorphism $\cl^+_{2,2}\simeq\cl_{2,1}\simeq\Om_{2,0}\simeq
\cl_{2,0}\oplus\cl_{2,0}$ we have for the ideal a following
reduction: $I_{2,2}\simeq
\cl^+_{2,2}\frac{1}{2}(1+\cE_{13})\frac{1}{2}(1+\cE_{24})\simeq\Om_{2,0}
\frac{1}{2}(1-i\Upsilon_{12})$, where $\Upsilon_i$ are matrix representations
of the units of $\Om_{2,0}$:
\[
\ar\Upsilon_1=\begin{pmatrix}
0 & 1\\
1 & 0
\end{pmatrix},\quad\Upsilon_2=\begin{pmatrix}
0 & -i\\
i & 0
\end{pmatrix}.
\]
Thus, in the case of the immersion $S^{2,0}\hookrightarrow M^{2,2}$ 
we have a restriction of the spinor field
$\Psi\in\cl_{2,2}\frac{1}{2}(1+\cE_{13})\frac{1}{2}
(1+\cE_{24})$ of the manifold $M^{2,2}$ onto a spinor field 
$\Psi_{|S^{2,0}}\in\Om_{2,0}\frac{1}{2}(1-i\Upsilon_{12})$ of the surface
$S^{2,0}$, where $\Psi_{|S^{2,0}}=(\epsilon^+_{20}\Psi,\epsilon^-_{20}\Psi)$,
and $\epsilon^\pm_{20}=\frac{1}{2}(1\pm i\zeta_1\zeta_2)$ 
are mutually orthogonal idempotents of the anti--quaternion 
$\sum^3_{i=0}\cl^i_{2,0}\zeta_i$. 

Analogously, in the case of the immersion $S^{0,2}\hookrightarrow M^{2,2}$ 
the general element of $\cl_{2,2}$ is represented by a quaternion
$\sum^3_{i=0}
\cl^i_{0,2}\zeta_i$, where $\zeta_1=\e_{134},\zeta_2=\e_{234},\,\zeta_1^2=
\zeta^2_2=-1$. At this point the quaternions $\cl^i_{0,2}$ are generated by
a set
$\{1,\e_3,\e_4,\e_{34}\}$. Further, in the case of the immersion $S^{1,1}
\hookrightarrow M^{2,2}$ we have for the general element of $\cl_{2,2}$
a following pseudo--quaternion $\sum^3_{i=3}\cl^i_{1,1}\zeta_i$, where
$\zeta_1=\e_{123},\zeta_2=\e_{134},\zeta^2_1=-\zeta^2_2=1,\,\cl^i_{1,1}
\simeq\{1,\e_2,\e_4,\e_{24}\}$.
Therefore, on the surfaces
$S^{0,2}\hookrightarrow M^{2,2}$ and $S^{1,1}\hookrightarrow M^{2,2}$
there exist restricted spinor fields
$\Psi_{|S^{0,2}}=(\epsilon^+_{02}\Psi,\epsilon^-_{02}\Psi)$ and $\Psi_{|S^{1,1}}=
(\epsilon^+_{11}\Psi,\epsilon^-_{11}\Psi)$, where $\Psi\in I_{2,2},\,
\epsilon^\pm_{02}=\frac{1}{2}(1\pm i\zeta_1\zeta_2),\epsilon^\pm_{11}=
\frac{1}{2}(1\pm\zeta_1\zeta_2)$.

Summarizing obtained above results we come to the following
\begin{theorem}\label{timmer}
Let $P=Q\times_{\spin(p,q)}\Delta_{p,q}$ be a spinor bundle of the
4--dimensional pseudo--riemannian manifold $M^{p,q}$ an let $\Psi\in\Gamma(P)$
be a smooth section (spinor field) of the bundle $P$. At this point the
components of $\Psi$ are elements of a minimal left ideal
$I_{p,q}$ of the Clifford algebra
$\cl_{p,q}$ of a tangent bundle of the manifold $M^{p,q}$. Then, a
restriction of the spinor bundle $P$ of $M^{p,q}$ onto a spinor
bundle $P_{|S^{p,q}}$ of a surface
$S^{p,q}$ ($p+q=2$) immersed into $M^{p,q}$ defined as follows:\\
1) At the immersions $S^{1,1}\hookrightarrow M^{1,3},\,S^{0,2}\hookrightarrow
M^{1,3}$ a spinor (mother) field 
$\Psi\in I_{1,3}=\cl_{1,3}\frac{1}{2}(1+\gamma_0)
$ on the Lorentzian manifold $M^{1,3}$ in virtue of the ideal reduction
$I_{1,3}=\cl^+_{1,3}\frac{1}{2}(1+\gamma_0)\simeq
\C_2\frac{1}{2}(1+i\sigma_{12})$ induces restricted spinor fields
$\Psi_{|S^{1,1}}=(\epsilon^+_{02}\Psi,\epsilon^-_{02}\Psi)$ and $\Psi_{|S^{0,2}}
=(\epsilon^+_{11}\Psi,\epsilon^-_{11}\Psi)$ on the surfaces $S^{1,1}$ and
$S^{0,2}$, where $\epsilon^\pm_{02}$ and $\epsilon^\pm_{11}$ are mutually
orthogonal idempotents respectively of the quaternion $\sum^3_{i=0}\cl^i_{1,1}
\zeta_i$ and pseudo--quaternion $\sum^3_{i=0}\cl^i_{0,2}\zeta_i$, by means of
which (in virtue of the decompositions $\cl_{1,3}\simeq\cl_{1,1}\otimes\cl_{0,2}$ 
and $\cl_{1,3}\simeq\cl_{0,2}\otimes\cl_{1,1}$) is represented a general
element of the spacetime algebra $\cl_{1,3}$. $\sigma_i$ and $\gamma_i$ are
respectively the Pauli and Dirac matrices, $\C_2$ is an algebra of
hyperbolic biquaternions.\\
2) At the immersions $S^{1,1}\hookrightarrow M^{3,1},\,S^{2,0}\hookrightarrow
M^{3,1}$ a spinor field $\Psi\in I_{3,1}=\cl_{3,1}\frac{1}{2}(1+\cE_{1})
\frac{1}{2}(1+\cE_{24})$ on the Majorana manifold $M^{3,1}$ by virtue of the
reduction
$I_{3,1}=\cl^+_{3,1}\frac{1}{2}(1+\cE_{1})\frac{1}{2}(1+\cE_{24})\simeq
\C_2\frac{1}{2}(1+i\sigma_{12})$ induces restricted spinor fields
$\Psi_{|S^{1,1}}=(\epsilon^+_{20}\Psi,\epsilon^-_{20}\Psi)$ and
$\Psi_{|S^{2,0}}=(\epsilon^+_{11}\Psi,\epsilon^-_{11}\Psi)$ on the surfaces
$S^{1,1}$ and $S^{2,0}$, where $\epsilon^\pm_{20}$ and $\epsilon^\pm_{11}$ are
mutually orthogonal idempotents respectively of the anti--quaternion
$\sum^3_{i=0}\cl^i_{1,1}\zeta_i$ and pseudo--quaternion 
$\sum^3_{i=0}\cl^i_{2,0}
\zeta_i$, by means of which (in virtue of the decompositions 
$\cl_{3,1}\simeq\cl_{1,1}
\otimes\cl_{2,0}$ and $\cl_{3,1}\simeq\cl_{2,0}\otimes\cl_{1,1}$)
is represented a general element of the Majorana algebra $\cl_{3,1}$.\\
3) At the immersions $S^{2,0}\hookrightarrow M^{4,0},\,S^{0,2}\hookrightarrow
M^{0,4}$ spinor fields $\Psi\in I_{4,0}=\cl_{4,0}\frac{1}{2}(1+\cE_{1})
,\,\Psi\in I_{0,4}=\cl_{0,4}\frac{1}{2}(1+\cE_{123})
$ of the manifolds $M^{4,0},\,M^{0,4}$ by virtue of the reductions
$I_{4,0}\simeq\cl^+_{4,0}\frac{1}{2}(1+\cE_{1})\simeq
\Om_{0,2}\frac{1}{2}(1+i\Upsilon_{12}),\,I_{0,4}=\cl^+_{0,4}\frac{1}{2}(1+
\cE_{123})\simeq\Om_{0,2}\frac{1}{2}(1+i\Upsilon_{12})$
induce restricted spinor fields $\Psi_{|S^{2,0}}=(\epsilon^+_{02}\Psi,
\epsilon^-_{02}\Psi),\,\Psi_{|S^{0,2}}=(\epsilon^+_{20}\Psi,\epsilon^-_{20}
\Psi)$ on the surfaces $S^{2,0},\,S^{0,2}$ respectively, where 
$\epsilon^\pm_{02}$ and $\epsilon^\pm_{20}$ are mutually orthogonal
idempotents of the quaternion $\sum^3_{i=0}\cl^i_{2,0}\zeta_i$ and
anti--quaternion
$\sum^3_{i=0}\cl^i_{0,2}\zeta_i$, by means of which (in virtue of the
decompositions
$\cl_{4,0}\simeq\cl_{2,0}\otimes\cl_{0,2}$ and $\cl_{0,4}\simeq\cl_{0,2}
\otimes\cl_{2,0}$) represented general elements of the algebras $\cl_{4,0}$ and
$\cl_{0,4}$. $\Upsilon_i$ are matrix representations of the units of an
algebra of elliptic biquaternions $\Om_{0,2}$.\\
4) At the immersions $S^{2,0}\hookrightarrow M^{2,2},\,S^{0,2}\hookrightarrow
M^{2,2},\,S^{1,1}\hookrightarrow M^{2,2}$ a spinor field $\Psi\in
I_{2,2}=\cl_{2,2}\frac{1}{2}(1+\cE_{13})\frac{1}{2}(1+\cE_{24})$
on the Kleinian manifold $M^{2,2}$ by virtue of the reduction
$I_{2,2}\simeq\cl^+_{2,2}\frac{1}{2}(1+
\cE_{13})\frac{1}{2}(1+\cE_{24})\simeq\Om_{2,0}\frac{1}{2}(1-i\Upsilon_{12})$
induces restricted spinor fields $\Psi_{|S^{2,0}}=(\epsilon^+_{20}\Psi,
\epsilon^-_{20}\Psi),\,\Psi_{|S^{0,2}}=(\epsilon^+_{02}\Psi,\epsilon^-_{02}
\Psi)$ and $\Psi_{|S^{1,1}}=(\epsilon^+_{11}\Psi,\epsilon^-_{11}\Psi)$
on the surfaces $S^{0,2},\,S^{2,0}$ and $S^{1,1}$, where $\epsilon^\pm_{20},\,
\epsilon^\pm_{02}$ and $\epsilon^\pm_{11}$ are mutually orthogonal
idempotents respectively of the anti--quaternion 
$\sum^3_{i=0}\cl^i_{2,0}\zeta_i$, 
quaternion $\sum^3_{i=0}\cl^i_{0,2}\zeta_i$ and pseudo--quaternion $\sum^3
_{i=0}\cl^i_{1,1}\zeta_i$, by means of which (in virtue of the
decompositions
$\cl_{2,2}\simeq\cl_{2,0}\otimes\cl_{2,0},\,\cl_{2,2}\simeq\cl_{0,2}\otimes
\cl_{0,2}$ and $\cl_{2,2}\simeq\cl_{1,1}\otimes\cl_{1,1}$) is represented
a general element of the algebra $\cl_{2,2}$.
\end{theorem}
\noindent{\bf Remark 1}.
In the case of the Lorentzian manifold the spinor fields
$\psi=\Psi_{|S^{r,s}}$ on the surfaces $S^{1,1}\hookrightarrow M^{1,3}$
and $s^{0,2}\hookrightarrow M^{1,3}$ in accordance with \cite{Cra85,Lou93}
may be expressed via bilinear covariants $\sigma,\,\bJ,\,\bS,\,\bK,\,
\omega$, i.e. $\psi\simeq Z\eta$, where $Z=\sigma+\bJ+i\bS-i\gamma_{0123}\bK+
\gamma_{0123}\omega$, $\eta$ is an arbitrary complex number, and
\begin{eqnarray}
\sigma&=&\psi^+\gamma_0\psi=4<\widetilde{\psi^\ast}\psi>_0,\nonumber\\
J_\mu&=&\psi^+\gamma_0\gamma_\mu\psi=4<\widetilde{\psi^\ast}\gamma_\mu\psi>_0,
\nonumber\\
S_{\mu\nu}&=&\psi^+\gamma_0i\gamma_{\mu\nu}\psi=4<\widetilde{\psi^\ast}
i\gamma_{\mu\nu}\psi>_0,\nonumber\\
K_\mu&=&\psi^+\gamma_0i\gamma_{0123}\gamma_\mu\psi=4<\widetilde{\psi^\ast}
i\gamma_{0123}\gamma_\mu\psi>_0,\nonumber\\
\omega&=&-\psi^+\gamma_0\gamma_{0123}\psi=-4<\widetilde{\psi^\ast}
\gamma_{0123}\psi>_0.\nonumber
\end{eqnarray}
At this point the bilinear covariants satisfy to {\it Fierz identities }
\begin{gather}
\bJ^2=\sigma^2+\omega^2,\quad\bK^2=-\bJ^2,\nonumber\\
\bJ\cdot\bK=0,\quad\bJ\mw\bK=-(\omega+\gamma_{0123}\sigma)\bS.\nonumber
\end{gather}
The spinor field $\psi$, whose $\sigma,\bJ,\bS,\bK,\omega$
satisfy to Fierz identities, recovered by its bilinear covariants with
an accuracy of the complex factor $\eta$. Moreover, both in the non--null
$(\sigma,\omega\neq 0)$ and null case
$(\sigma,\omega=0)$ the spinor $\psi$ is defined by bilinear covariant
$Z$ ($\psi=(1/4N)e^{-i\alpha}Z\eta)$, where $N=\sqrt{<Z\eta>_0}$) which
in its turn is defined by the spinor $\psi$ as follows:
$Z=4\psi\widetilde{\psi^\ast}=4\psi\psi^+\gamma_0$. Thus, we have a
so--called {\it boomerang} \cite{Lou93}. All bilinear covariants are real
and have important meaning in the Dirac theory of electron. 
In perspective, it is of interest to consider the analogous bilinear
covariants and boomerangs for the spinor fields on the surfaces immersed
into 4d manifolds with signatures different from the signature of the
Lorentzian manifold.\\[0.3cm]
\noindent{\bf Remark 2}. In more general case of non--orientable manifolds
we come to a group $\pin(p,q)$ which is a double covering of the structure
group $O(p,q)$ of the manifold $M^{p,q}$. In accordance with \cite{Dab88,
BD89} there exist eight double coverings of the orthogonal group $O(p,q)$:
\[
\rho^{a,b,c}:\,\pin^{a,b,c}(p,q)\simeq\frac{(\spin_0(p,q)\odot C^{a,b,c})}
{\Z_2}\longrightarrow O(p,q),
\]
where $C^{a,b,c}\in\{\Z_2\otimes\Z_2\otimes\Z_2,\,\Z_2\otimes\Z_4,\,Q_4,\,
D_4\}$ is a double covering of a discrete group of the space $\R^{p,q}$,
$a,b,c\in\{+,-\}$. In connection with this it is of interest to define
$\pin^{a,b,c}$--fields (generalization of ordinary spinor fields) on surfaces
immersed into the non--orientable manifolds. The classification of these
fields may be easily defined with the usage of recently established
relation between signatures of the spaces $\R^{p,q}$ and finite groups of
fundamental automorphisms of the Clifford algebras (see 
\cite[theorem 10]{Var99c}).\section{The Dirac Operator on the Surfaces Immersed into 4D Manifolds}
Let us consider now the Dirac operator on the surfaces 
$S^{r,s}\hookrightarrow M^{p,q}\;(r+s=2,\,p+q=4)$. First of all, let recall
some basic facts about a theory of the Dirac operator on a spin manifold
\cite{BFGK}. Let $(M^{p,q},g)$ be a pseudo--riemannian spin manifold and
let $S=Q\times_{\spin(p,q)}\Delta_{p,q}$ be a spinor bundle of the
manifold $(M^{p,q},g)$.
{\it The Dirac operator} on the manifold $(M^{p,q},g)$ is a first order
self--adjoint elliptic differential operator defined by an expression
\begin{equation}\label{DefDir}
D:\;\Gamma(S)\stackrel{\nabla^S}{\longrightarrow}\Gamma(TM\otimes S)
\stackrel{\mu}{\longrightarrow}\Gamma(S).
\end{equation}
where $\mu$ is a so--called {\it Clifford multiplication}:
\begin{eqnarray}
\mu:&&\R^n\otimes\Delta_n\longrightarrow\Delta_n\nonumber\\
&&x\otimes u\longmapsto\mu(x\otimes u)=x\cdot u=\nonumber\\
&&\begin{cases}
\gamma_{p,q}(x)u, & \text{if $p-q=0,2,4,6\pmod{8}$};\\
\text{proj}_j\cdot\gamma_{p,q}(x)u, & \text{if $p-q=1,3,5,7\pmod{8}$}.
\end{cases}\label{clif}
\end{eqnarray}
where $\gamma:\,\cl_{p,q}\rightarrow\End_{\K}(I_{p,q})$ is a spinor
representation of the Clifford algebra $\cl_{p,q}$. For the spinor bundles
the Clifford multiplication is defined as follows
\begin{eqnarray}
\mu:&&TM\otimes S\longrightarrow S\nonumber\\
&&x\otimes\varphi\longmapsto x\cdot\varphi\nonumber
\end{eqnarray}
In the case of even dimensions $(p-q\equiv\pmod{0,2,4,6})$ the Clifford
multiplication $\mu$ exchanges the positive and negative parts of
the bundle $S$. The Clifford multiplication may be also defined for the
$k$--forms. The action of the $k$--form
$\omega\in\Omega^k(M)$ on the spinor bundle is defined by a following
local formula
\[
\omega\cdot\varphi=\sum_{1\leq i_1<\ldots<i_k\leq n}\omega(s_{i_1},\ldots
,s_{i_k})s_{i_1}\cdot\ldots\cdot s_{i_k}\cdot\varphi,
\]
where $(s_1,\ldots,s_n)$ is a local orthonormal basis of the manifold
$(M^n,g)$, $\varphi\in\Gamma(S)$, $n=p+q$.

Further, 
$\nabla^S:\,\Gamma(S)\longrightarrow\Gamma(TM
\otimes S)$ in (\ref{DefDir}) is {\it a spinor derivative}, which locally
is given by an expression
\begin{equation}\label{loc}
\nabla^S_X\varphi=X(\varphi)+\frac{1}{2}\sum_{1\leq k<l\leq n}\omega_{kl}
(x)s_k\cdot s_l\cdot\varphi,
\end{equation}
where $\omega_{kl}=g(\nabla^Ms_k,s_l)$ are the connection forms of the
Levi--Civita connection $\nabla^M$ on $(M^{p,q},g)$ with respect to a
local basis $(s_1,\ldots,s_n)$, $X$ is a vector field.

It immediately follows that locally the Dirac operator may be written in
the form
\begin{equation}\label{Dirop}
D=\sum^n_{k=1}s_k\cdot\nabla^S_{s_k}.
\end{equation}

Since the all 4--dimensional manifolds are quaternionic manifolds, then
in the each point of such a manifold the Clifford algebra of the tangent
space is isomorphic to a quaternionic algebra, i.e. a Clifford bundle of the
manifold in this case may be represented in terms of the  quaternionic
algebras.
Indeed, in the case of even dimension, the volume
element $\omega=\e_{12\ldots n}$ is not belong to a center of the algebra
$\cl_{n}$. However, when $i\leq 2m$ we have
\begin{eqnarray}
\e_{12\ldots 2m 2m+k}\e_{i}&=&(-1)^{2m+1-i}\sigma(i-l)
\e_{12\ldots i-1 i+1\ldots 2m 2m+k},\nonumber \\
\e_{i}\e_{12\ldots 2m 2m+k}&=&(-1)^{i-1}\sigma(i-l)
\e_{12\ldots i-1 i+1\ldots 2m 2m+k},\nonumber
\end{eqnarray}
where $\sigma(n)$ are the functions of the form (\ref{e2}).
Therefore, the commutativity condition of the elements 
$\e_{12\ldots 2m 2m+k}$
and $\e_{i}$ is $2m+1-i\equiv i-1\pmod{2}$. Thus, the elements 
$\e_{12\ldots 2m 2m+1}$ and $\e_{12\ldots 2m 2m+2}$ are commute with all basis
elements $\e_{i}$ whose indexes are not exceed $2m$. Therefore, a transition
from $\cl_{2m}$ to $\cl_{2m+2}$  
may be represented as transition
from the real coordinates in the algebra
$\cl_{2m}$
to quaternionic coordinates of the form $a+b\zeta_1+c\zeta_2+d\zeta_1\zeta_2$, 
where $\zeta_1$ and $\zeta_2$
are additional basis elements $\e_{12\ldots 2m 2m+1}$ and 
$\e_{12\ldots 2m 2m+2}$. The elements $\e_{i_{1}i_{2}\ldots i_{k}}\zeta_1$
are contain index $2m+1$ and not contain index $2m+2$, and the elements
$\e_{i_{1}i_{2}\ldots i_{k}}\zeta_2$ are contain index $2m+2$ and not contain
index $2m+1$. Respectively, the elements $\e_{i_{1}i_{2}\ldots i_{k}}
\zeta_1\zeta_2$ are contain both indexes $2m+1$ and $2m+2$. 
Therefore, the algebras
$\cl_{p,q+2},\,\cl_{p+2,q}$ and $\cl_{p+1,q+1}\;(p-q\equiv 0,2,4,6\pmod{8})$
are isomorphic respectively to quaternionic, anti--quaternionic and
pseudo--quaternionic algebras, i.e. a general element of these algebras
can be represented in the form
\begin{equation}\label{elem}
\cl^0_{p,q}+\cl^1_{p,q}\zeta_1+\cl^2_{p,q}\zeta_2+\cl^3_{p,q}\zeta_1\zeta_2,
\end{equation}
where $\zeta_1=\e_{12\ldots 2m 2m+1},\,\zeta_2=\e_{12\ldots 2m 2m+2}$.
Respectively, in dependence on the squares of the units $\zeta_1,\zeta_2$
the expression (\ref{elem}) is the quaternion $(\zeta^2_1=\zeta^2_2=-1)$,
anti--quaternion $(\zeta^2_1=\zeta^2_2=1)$ and pseudo--quaternion
$(\zeta^2_1=-\zeta^2_2=1)$. In other words, according to the theorem
\ref{tKar} we have for the quaternionic algebras the following
decompositions:
\begin{eqnarray}
\cl_{p,q+2}&\simeq&\cl_{0,2}\otimes\cl_{q,p},\nonumber\\
\cl_{p+2,q}&\simeq&\cl_{2,0}\otimes\cl_{q,p},\label{gendec}\\
\cl_{p+1,q+1}&\simeq&\cl_{1,1}\otimes\cl_{p,q}.\nonumber
\end{eqnarray}
These decompositions are natural generalizations of the decompositions
considered above in the section 3.

Let $\spin(2m)\subset\cl^\star_{2m}$ be a spinor group
and let $\Delta_{2m}=\Delta^+_{2m}\oplus
\Delta^-_{2m}$ be a representation of the group $\spin(2m)$, $\cl^\star_{2m}$
is a set of all invertible elements of the algebra $\cl_{2m}$.
\begin{lem}\label{lrest1}
The restriction of $\Delta_{2m+2}$ to $\spin(2m)$ is isomorphic to the
$\spin(2m)$--representation $\Delta_{2m}$, where an action of 
$\zeta_1\zeta_2$ on $\Delta_{2m}=\Delta^+_{2m}\oplus\Delta^-_{2m}$
is defined by an expression
\[
\zeta_1\zeta_2\cdot(u^+\oplus u^-)=(-1)^m\varepsilon u^+-(-1)^m\varepsilon
u^-.
\]
Here $\varepsilon=i$ if $\zeta^2_1=\zeta^2_2=\pm 1$ and
$\varepsilon=1$ if $\zeta^2_1=-\zeta^2_2=1$.
\end{lem}
\begin{proof} The spinor group $\spin(2m+2)$ is completely defined in terms
of the algebra $\cl_{2m+2}$:
\begin{equation}\label{cl1}
\spin(2m+2)=\left\{s\in\Gamma^+_{2m+2}|N(s)=\pm 1\right\},
\end{equation}
where $s\in\cl^\star_{2m+2},\,N:\cl_{2m+2}\rightarrow\cl_{2m+2},\,N(x)=x
\widetilde{x};\,\Gamma^+_{2m+2}=\Gamma_{2m+2}\cap\cl^+_{2m+2}$ is a special
Clifford--Lipschitz group, and
\begin{equation}\label{cl2}
\Gamma_{2m+2}=\left\{s\in\cl_{2m+2}|\forall x\in\R^{2m+2},sx\hat{s}^{-1}\in
\R^{2m+2}\right\}.
\end{equation}
Let $\rho^{2m+2}:\,\cl_{2m+2}\rightarrow\End E$ is a representation of the
algebra $\cl_{2m+2}$ in a vector space $E$. The representation $\rho^{2m+2}$
induces via (\ref{cl2}) a representation of the group $\pin(2m+2)=\left\{
s\in\Gamma_{2m+2}|N(s)=\pm 1\right\}$, and also via (\ref{cl1}) a
representation $\Delta_{2m+2}$ of the group $\spin(2m+2)$. Further, in 
virtue of the decomposition (\ref{elem}) an inverse transition
$\cl_{2m+2}\rightarrow\cl_{2m}$ induces a transition $\spin(2m+2)\rightarrow
\spin(2m)$, and $\rho^{2m+2}\rightarrow\rho^{2m}$ induces a restriction
$\Delta_{2m+2}\rightarrow\Delta_{2m}$
by means of mutually orthogonal idempotents
(projection operators) $\epsilon^{\pm}=\frac{1}{2}
(1\pm\varepsilon\zeta_1\zeta_2)$. At this point $\zeta_1\zeta_2=\e_{2m+1\, 2m+2}
\mapsto\cE_{2m+1\, 2m+2}$ commutes with $u^+\in\Delta^+_{2m}$ and 
anticommutes with $u^-\in\Delta^-_{2m}$.\end{proof}

Further, let $M$ be an $(2m+2)$--dimensional pseudo--riemannian manifold and
let $F$ be an $2m$--dimensional submanifold immersed into $M$,
$F\hookrightarrow M$.
We suppose that the both manifolds endowed with some spinor structure.
Let $N$ be a normal bundle of the manifold
$F\hookrightarrow M$, then in accordance with \cite{Mil65} a sum 
of the spinor structures on the tangent bundle and on the normal bundle
of the manifold $F$ coincides with the spinor structure on the tangent
bundle of the manifold $M$ restricted to $F$. Let $\nabla^F$ and $\nabla^M$ be
Levi--Civita connections on the manifolds $F$ and $M$, respectively. Let
$\nabla^N$ be a normal connection on the bundle $N$. Denote the second
fundamental form of the submanifold $F^{2m}\hookrightarrow M^{2m+2}$
as $II$. Further, let $\zeta_1$ and $\zeta_2$ be unit normal vector fields
on $F^{2m}\hookrightarrow M^{2m+2}$ and let
$S_F=Q_F\times_{\spin(2m)}\Delta_{2m}$ be a spinor bundle of the
submanifold $F^{2m}$.
\begin{lem}\label{lrest2}
If $n=2m+2$, then a restriction of the spinor bundle $S$ of the manifold
$(M^n,g)$ onto submanifold $F^{n-2}$ is isomorphic to the bundle $S_F$, where
$\zeta_1\zeta_2$ acts on $S_F$ as follows
\begin{equation}\label{act}
\zeta_1\zeta_2\cdot(\psi^+\oplus\psi^-)=(-1)^m\varepsilon\psi^+-
(-1)^m\varepsilon\psi^-,
\end{equation}
and the spinor derivative of $\psi\in\Gamma(S)$ equals
\begin{equation}\label{sder}
\nabla^S_X\psi=\left(\nabla^{S_F}_X\otimes\Id+
\Id\otimes\nabla^{S_N}_X\right)\psi+\frac{1}{2}
\sum^{2m}_{1\leq i_1<\ldots<i_k\leq 2m}<II(X,X_i),X_i>\cdot\zeta_1\zeta_2
\cdot\psi
\end{equation}
for all $X\in T_xF$.
\end{lem}
\begin{proof} The expression (\ref{act}) immediately follows from the lemma
\ref{lrest1}. Further, follows to \cite{Bar97c} we see that for some
point $x\in F$ and a vector field
$X\in T_xF$ the Gauss formula with respect to a decomposition
$T_xM=T_xF\oplus N_x$ gives
\[
\nabla^M_X=\begin{pmatrix}
\nabla^F_X & -II(X,\cdot)^\ast\\
II(X,\cdot) & \nabla^N_X
\end{pmatrix},
\]
or
\begin{equation}\label{FGauss}
\nabla^M_X-\left(\nabla^F_X\oplus\nabla^N_X
\right)=\begin{pmatrix}
0 & -II(X,\cdot)^\ast\\
II(X,\cdot) & 0
\end{pmatrix}.
\end{equation}

Let $X_1,\ldots,X_{2m}$ be a local orthonormal tangent frame of the
submanifold $F^{2m}$ at the point $x$ and let $Y_1,Y_2$ be a local
orthonormal frame on the normal bundle $N$ at $x$. Then 
$h:=(X_1,\ldots,X_{2m},\,
Y_1,Y_2)$ is a local section of the tangent bundle $P\times_{SO(2m+2)}$
of the manifold $M^{2m+2}$ restricted to $F^{2m}$. Now we can to write
(\ref{FGauss}) in the matrix form:
\begin{equation}\label{FGauss2}
\nabla^M_X-\left(\nabla^F_X\oplus\nabla^N_X
\right)=\begin{pmatrix}
0 & (-<II(X,X_i),Y_j>)_{j,i}\\
(<II(X,X_i),Y_j>)_{i,j} & 0
\end{pmatrix}.
\end{equation}
Further, let $\bom^F,\,\bom^N$ and $\bom^M$ be respectively connection 1--forms
for $\nabla^F,\,\nabla^N$ and $\nabla^M$ lifted to 
$\spin(2m),\,\spin(2)$ and $\spin(2m+2)$.
If $\Theta:\;\spin(2m+2)\rightarrow SO(2m+2)$ is an usual double covering,
then (\ref{FGauss2}) can be written as follows
\begin{multline}\label{FG3}
\Theta_\ast(\bom^M(dh\cdot X)-(\bom^F\oplus\bom^N)(dh\cdot X)=\\
\begin{pmatrix}
0 & (-<II(X,X_i),Y_j>)_{j,i}\\
(<II(X,X_i),Y_j>)_{i,j} & 0
\end{pmatrix}.
\end{multline}
Using the standard formula \cite[c.42]{LM89} for $\Theta_\ast$ we obtain from
(\ref{FG3}) 
\begin{equation}\label{FG4}
\bom^M(dh\cdot X)-(\bom^F\oplus\bom^N)(dh\cdot X)=\frac{1}{2}
\sum^{2m}_{i=1}\sum^2_{j=1}<II(X,X_i),Y_j>\cdot e_i\cdot f_j,
\end{equation}
where $e_1,\ldots,e_{2m}$ is a standard basis of the space $\R^{2m}$,
$f_1,f_2$ is a standard basis of the space $\R^2$. 

Let $S_M=Q\times_{\spin(2m+2)}\Delta_{2m+2}$ be a spinor bundle of the
manifold $M^{2m+2}$. In virtue of the decomposition 
$\cl_{2m+2}\simeq\cl_2\otimes
\cl_{2m}$ we have $S_{M|F}=S_F\otimes S_N$, where $S_F=Q_F\times_{\spin
(2m)}\Delta_{2m},\,S_N=Q_N\times_{\spin(2)}\Delta_2$. Let
$\nabla^{S_M},\,\nabla^{S_F}$ and $\nabla^{S_N}$ be Levi--Civita connections
on the bundles $S_M,\,S_F$ and $S_N$, respectively. Then
\[
\nabla^{S_F\otimes S_N}:=\nabla^{S_F}\otimes\Id+\Id\otimes\nabla^{S_N}
\]
is a Levi--Civita product connection on $S_F\otimes S_N$. At this point
the equation (\ref{FG4}) takes a form
\begin{equation}\label{FG5}
\nabla^S_X-\left(\nabla^{S_F}_X\otimes\Id+\Id\otimes\nabla^{S_N}_X\right)=
\frac{1}{2}\sum^{2m}_{i=1}\sum^2_{j=1}<II(X,X_i),Y_j>\mu(X_i\cdot Y_j),
\end{equation}
where $\mu(X_i\cdot Y_j)$ is the Clifford multiplication defined by
(\ref{clif}). Whence in accordance with the definition of the spinor
derivative (\ref{loc}) and identifications
$\zeta_1\leftrightarrow Y_1,\zeta_2\leftrightarrow Y_2$
follows the formula (\ref{sder}).
\end{proof}

Before defining the Dirac operator (by the formula (\ref{Dirop})), 
corresponding to the spinor derivative (\ref{FG5}) it is necessary to consider
the following two operators
\begin{equation}\label{Bar1}
\widetilde{D}=\sum^{2m}_{j=1}X_j\cdot\nabla^{S_F\otimes S_N}_{X_j}
\end{equation}
and
\begin{equation}\label{Bar2}
\hat{D}=\sum^{2m}_{j=1}X_j\cdot\nabla^S_{X_j}.
\end{equation}
It is easy to see that both operators act on the sections of the
bundle $S_M$. 
Let $H=\frac{1}{2m}\sum^{2m}_{j=1}II(X_j,X_j)$ be the mean curvature vector
field of the submanifold $F^{2m}\hookrightarrow M^{2m+2}$. Further,
using (\ref{FG5}) we obtain
\begin{eqnarray}
\hat{D}-\widetilde{D}&=&\frac{1}{2}\sum^{2m}_{i,j=1}X_j\cdot
<II(X_j,X_i),X_i>\cdot\zeta_1\zeta_2\nonumber\\
&=&\frac{1}{2}\sum^{2m}_{i,j=1}\mu(X_j\cdot X_i)<II(X_j,X_i)>\cdot
\zeta_1\zeta_2.\nonumber
\end{eqnarray}
By virtue of antisymmetry, the products $X_j\cdot X_i$ with $i\neq j$ vanish,
at this point the form $II(X_j,X_i)$ is symmetric. Since
$X\cdot Y+Y\cdot X=2g(X,Y)\Id_{S_M}$, then
\begin{equation}\label{Bar3}
\hat{D}-\widetilde{D}=\frac{1}{2}\sum^{2m}_{i=1}g_{ii}H\cdot\zeta_1\zeta_2.
\end{equation}
\begin{theorem}\label{tDir}
Let $\Psi\in I_{p,q}$ be a real Killing spinor field on the 4--dimensional
pseudo--riemannian manifold $M^{p,q}$ and let $\Psi_{|S^{r,s}}=\psi=\psi^+
\oplus\psi^-$ be a restricted spinor field on the surface $S^{r,s}$
immersed into the manifold $M^{p,q}$. Then a Dirac operator of the surface
$S^{r,s}\hookrightarrow M^{p,q}$ defined as follows
\begin{eqnarray}
D(\psi^+)&=&(\alpha-\frac{1}{2}\varepsilon\beta H)\psi^-,\nonumber\\
D(\psi^-)&=&(\alpha+\frac{1}{2}\varepsilon\beta H)\psi^+,\nonumber
\end{eqnarray}
where $\varepsilon=1$ for the immersions $S^{0,2}\hookrightarrow M^{1,3},\,
S^{2,0}\hookrightarrow M^{3,1},\,S^{1,1}\hookrightarrow M^{2,2}$, and
$\varepsilon=i$ for the immersions $S^{2,0}\hookrightarrow M^{4,0},\,
S^{0,2}\hookrightarrow M^{0,4},\,S^{1,1}\hookrightarrow M^{1,3},\,
S^{1,1}\hookrightarrow M^{3,1},\,S^{2,0}\hookrightarrow M^{2,2},\,
S^{0,2}\hookrightarrow M^{2,2}$. $H$ is a mean curvature of the surface,
$\alpha=\lambda_1g_{11}+\lambda_2g_{22},\,\beta=g_{11}+g_{22}$.
\end{theorem}
\begin{proof} In the case of the Lorentzian manifold $M^{1,3}$ we have
the following immersions $S^{0,2}\hookrightarrow M^{1,3}$ and $S^{1,1}
\hookrightarrow M^{1,3}$. At this point in accordance with the
theorem \ref{timmer} on the surfaces
$S^{0,2}$ and $S^{1,1}$ there exist the spinor fields
$\Psi_{|S^{0,2}}=(\psi^+,\psi^-)=(\epsilon^+_{11}\psi,\epsilon^-_{11}\psi)$
and $\Psi_{|S^{1,1}}=(\psi^+,\psi^-)=(\epsilon^+_{02}\psi,\epsilon^-_{02}\psi)$,
respectively.
Let us find a Dirac operator on the surface $S^{0,2}\hookrightarrow M^{1,3}$.
First of all, in accordance with the definition (\ref{Dirop}) and the
formulae (\ref{Bar1}), (\ref{Bar2}) let consider the following two operators:
\[
\widetilde{D}=X_1\cdot\nabla^{S_{S^{0,2}}\otimes S_N}_{X_1}+X_2\cdot
\nabla^{S_{S^{0,2}}\otimes S_N}_{X_2}
\]
and
\[
\hat{D}=X_1\cdot\nabla^{S_{M^{1,3}}}_{X_1}+X_2\cdot\nabla^{S_{M^{1,3}}}_{X_2}.
\]
The both operators act on the section of the spinor bundle 
$S=Q\times_{\spin(1,3)}
\Delta_{1,3}$. Using (\ref{sder}) and (\ref{Bar3}) we obtain
\begin{equation}\label{Dir1}
\hat{D}-\widetilde{D}=\frac{1}{2}\sum^2_{i,j=1}X_j\cdot X_i\cdot II(X_j,X_i)
\cdot\zeta_1\zeta_2=\frac{1}{2}\beta H\cdot\zeta_1\zeta_2,
\end{equation}
where $\beta=g_{11}+g_{22}$.
Further, let $\Psi\in I_{1,3}$ be the spinor field on the Lorentzian
manifold $M^{1,3}$, then from (\ref{Dir1}) follows
\begin{equation}\label{Dir2}
X_1\cdot\nabla^{S_{S^{0,2}}\otimes S_N}_{X_1}(\Psi)+X_2\cdot
\nabla^{S_{S^{0,2}}\otimes S_N}_{X_2}(\Psi)=D(\psi)-\frac{1}{2}
\beta H\cdot\zeta_1\zeta_2\cdot\psi,
\end{equation}
where $D(\psi)$ is a Dirac operator of the surface $S^{0,2}\hookrightarrow
M^{1,3}$ defined on the restriction $\psi=\Psi_{|S^{0,2}}$. We suppose now
that the spinor field $\Psi$ on the manifold $M^{1,3}$ is a real Killing
spinor, i.e. there exists such a number $\lambda\in\R$ that for any vector
field $X\in T_xM^{1,3}$ the derivative of $\Psi$ in the direction
$X$ equals
\[
\nabla^{S_{S^{0,2}}\otimes S_N}_X(\Psi)=\lambda\cdot X\cdot\Psi.
\]
Therefore, from (\ref{Dir2}) for the restriction $\psi=\Psi_{|S^{0,2}}$
we have
\[
D(\psi)=\alpha\psi+\frac{1}{2}\beta H\cdot\zeta_1\zeta_2\cdot\psi.
\]
where $\alpha=\lambda_1g_{11}+\lambda_2g_{22}$. Since 
$\psi=\psi^+\oplus\psi^-$, then in virtue of the relation (\ref{act}) of the
lemma \ref{lrest2} from the last equation we obtain (recalling that
$\zeta_1$ and $\zeta_2$ are the units of the pseudo--quaternion)
\begin{eqnarray}
D(\psi^+)&=&(\alpha-\frac{1}{2}\beta H)\psi^-,\nonumber\\
D(\psi^-)&=&(\alpha+\frac{1}{2}\beta H)\psi^+.\label{Dir3}
\end{eqnarray}
Further, for the immersion of the time--like surface $S^{1,1}\hookrightarrow
M^{1,3}$ the analogous calculations give (at this point $\zeta_1$ and
$\zeta_2$ are the quaternion units)
\begin{eqnarray}
D(\psi^+)&=&(\alpha-\frac{i}{2}\beta H)\psi^-,\nonumber\\
D(\psi^-)&=&(\alpha+\frac{i}{2}\beta H)\psi^+.\nonumber
\end{eqnarray}
\end{proof} 
The theorem \ref{tDir} has three important particular cases
\begin{cor}
Let $S^{r,s}\hookrightarrow M^{p,q}$ be a minimal surface, then a Dirac
operator of $S^{r,s}$ has a form
\[
D(\psi)=\alpha\psi,
\]
where $\psi=\psi^+\oplus\psi^-=\Psi_{|S^{r,s}}$ is an eigenspinor on
the surface $S^{r,s}$.
\end{cor}
\begin{cor}
Let $\Psi$ be a parallel spinor field $(\lambda_1=\lambda_2=0)$ on the
manifold $M^{p,q}$ and let $\psi=\Psi_{|S^{r,s}}$ be its restriction
on the surface $S^{r,s}\hookrightarrow M^{p,q}$, then a Dirac operator
of $S^{r,s}$ takes a form
\begin{eqnarray}
D(\psi^+)&=&-\frac{1}{2}\varepsilon\beta H\psi^-,\nonumber\\
D(\psi^-)&=&\frac{1}{2}\varepsilon\beta H\psi^+.\nonumber
\end{eqnarray}
\end{cor}
On the other hand, when $\varepsilon=i$ and $\beta=2,\,\lambda_1=\lambda_2=0$
(parallel spinor field) a Dirac operator on the surface may be written in
more compact form. Let consider a following spinor field
\begin{multline}
\psi^\circ=\psi^+-i\psi^-=\frac{1}{2}(\psi+i\cdot\zeta_1\zeta_2\cdot\psi)-
\frac{i}{2}(\psi-i\cdot\zeta_1\zeta_2\cdot\psi)=\\
\frac{1}{2}(1-i)\psi+
\frac{1}{2}(-1+i)\cdot\zeta_1\zeta_2\cdot\psi,\nonumber
\end{multline}
where $\psi\in\Om_{2,0}\frac{1}{2}(1-i\Upsilon_{12})$ for $S^{2,0}
\hookrightarrow M^{2,2}$ and $\psi\in\Om_{0,2}\frac{1}{2}(1+i\Upsilon_{12})$
for $S^{2,0}\hookrightarrow M^{4,0}$. Then
\[
D(\psi^\circ)=H\psi^\circ.
\]
Analogously, when $\varepsilon=i$ and $\beta=-2,\,\lambda_1=\lambda_2=0$
we have for the immersions $S^{0,2}\hookrightarrow M^{2,2},\,S^{0,2}
\hookrightarrow M^{0,4}$ a following spinor field
\[
\psi^\bullet=\psi^++i\psi^-=\frac{1}{2}(\psi+i\cdot\zeta_1\zeta_2\cdot\psi)+
\frac{i}{2}(\psi-i\cdot\zeta_1\zeta_2\cdot\psi)=\frac{1}{2}(1+i)(\psi+
\zeta_1\zeta_2\cdot\psi)
\]
and
\[
D(\psi^\bullet)=H\psi^\bullet.
\] 
\begin{cor}\label{c3}
If $\lambda_1=\lambda_2$, then in the case of the immersions of time--like
surfaces $S^{1,1}\hookrightarrow M^{p,q}$ a Dirac operator of $S^{1,1}$ is
homogeneous,
\[
D(\psi)=0.
\]
\end{cor}
\noindent{\bf Example}.
Let $S^{1,1}\hookrightarrow M^{1,3}$ be an immersion of the time--like
surface into the Lorentzian manifold and let $\lambda_1=\lambda_2$, then
$\alpha=\beta=0$ (corollary \ref{c3}) and
\[
D(\psi)=0.
\]
At this point it is easy to trace a relation with a so--called 
{\it optical geometry}
\cite{Rob61,Pen83,RT86,Nur96}. Indeed, let $\phi\in\cl^+_{1,3}
\simeq\cl_{3,0}\simeq\C_2$ be a Dirac--Hestenes spinor field an let
${\bf F}={\bf E}+i{\bf B}=\partial\mw A\in\cl_{3,0}\simeq\C_2$ 
be an electromagnetic (in general case, non--null) field,
where $\partial=\partial^0+
\partial^1\e_1+\partial^2\e_2+\partial_3\e_3$ and 
$A=A^0+A^1\e_1+A^2\e_2+A^3\e_3$ are partial derivative and vector--potential,
respectively. By virtue of $\sigma:\,\cl_{3,0}\rightarrow\End_{\C}(I_{3,0})$
the element $\Phi=F_1\e_1+F_2\e_2+F_3\e_{12}\in\C_2$ 
in the spinor representation is defined by a following symmetric matrix
\[
\Phi=\ar\begin{pmatrix} F_1+iF_2 & iF_3\\
iF_3 & F_1-iF_2\end{pmatrix}.
\]
The determinant $\det\Phi=F^2_1+F^2_2+F^2_3$ vanishes if, and only if, 
the electromagnetic field is null, i.e. when 
${\bf E}\cdot{\bf B}=0$ and ${\bf E}^2={\bf B}^2$. Null electromagnetic
fields play a key role in the theory of shear free congruences of null
geodesics in the Lorentzian manifold and give rise to the optical geometry
and a Cauchy--Riemann structure on the space of null geodesics.
Expressing the Dirac--Hestenes spinor field
$\phi\in\cl^+_{1,3}$ via the null electromagnetic field,
$\phi_1=\alpha+iB_3,\,\phi_2=-B_2+iB_1,\,\phi_3=E_3+i\lambda,
\phi_4=E_1+iE_2$ (see \cite{Par92}), we find that in the case of the
immersion of the time--like surface (or, light cone) $S^{1,1}$ into the
Lorentzian manifold $M^{1,3}$ a restriction of the spinor field 
$\Psi\in M^{1,3}$
onto a spinor field $\psi=\Psi_{|S^{1,1}}$ of the surface
$S^{1,1}\hookrightarrow M^{1,3}$ in accordance with (\ref{restrict2}) 
is expressed via the null electromagnetic field, and a Dirac operator
of $S^{1,1}$ in accordance with corollary \ref{c3} is homogeneous. Such a
form of the Dirac operator corresponds to massless physical fields, which
describe, as known, such particles as photon and neitrino.
\section{Local Spinor Representations of Surfaces in 4D Pseudo--Euclidean
spaces}
Let us consider now local (spinor) representations of surfaces conformally
immersed into 4--dimensional pseudo--euclidean spaces. As known, these
representations are defined by so--called Gauss map (GGM) 
\cite{HO80,HO83,HO85} and generalized Weierstrass representation (GWR) \cite{KL98a}.
The our main goal in this section is a realization of GGM and GWR in terms
of the spinor fields introduced above in the section 4.

Let $S^{r,s}\hookrightarrow M^{p,q}\,(r+s=2,p+q=4)$ be a surface endowed
with some spinor structure and let $S_0$ be a connected Riemann surface with
a local complex coordinate $z$. Let $P\times_G$
be a principal bundle of $S_0$ with the structure group $G$ ($G$ is a group
of fractional linear transformations). Then the {\it spinor representation
of a surface} $S^{r,s}$ in $M^{p,q}$ (or, locally, in $\R^{p,q}$) is given by  
the following diagramm
\[
\dgARROWLENGTH=0.9em
\begin{diagram}
\node{Q\times_{\widetilde{G}}}\arrow{s,l}{\mu}\arrow{e,t}{\chi}
\node{Q\times_{\spin(r,s)}}\arrow{s,r}{f}\\
\node{P\times_{G}}\arrow{s}\arrow{e,t}{\omega}\node{P\times_{SO(r,s)}}
\arrow{s}\\
\node{S_0}\arrow{e,t}{g}\node{S^{r,s}}
\end{diagram}
\]
Here $g:\,S_0\rightarrow G_{2,4}\simeq Q_2\simeq\C P^1\times\C P^1\simeq
S^2\times S^2$ is a generalized Gauss map. 

We start a consideration with the immersion of the space--like surface 
$S^{2,0}$ into the manifold
$M^{4,0}$. Locally, in virtue of an isomorphism 
$T_xM^{4,0}\simeq\R^{4,0}$ we have an immersion
$S^{2,0}\hookrightarrow\R^{4,0}$. The Grassmannian of oriented two--planes
in $\R^{4,0}$ may be identified with a quadric
$Q_2\subset\C P^3$, where $Q_2\simeq S^2\times S^2$, $S^2$ is a standard
sphere of radius $1/\sqrt{2}$. In this case, according to \cite{HO83}
generalized Gauss map $g:\,S_0\rightarrow G_{2,4}\simeq Q_2$
can be parametrized in terms of two complex functions
$f_1$ and $f_2$ as follows
\begin{equation}\label{Gmap}
\Phi(z)=(1+f_1f_2,\,i(1-f_1f_2),\,f_1-f_2,\,-i(f_1+f_2)),
\end{equation}
It is easy to see that $\sum^4_{k=1}\varphi^2_k=0$. The functions $f_1$
and $f_2$ are related by the formulae \cite{HO85}:
\begin{gather}
\left|\frac{f_{1\bar{z}}}{1+|f_1|^2}\right|=\left|\frac{f_{2\bar{z}}}
{1+|f_2|^2}\right|,\nonumber\\
\im\left\{\left(\frac{f_{1z\bar{z}}}{f_{1\bar{z}}}-\frac{2\bar{f}_1f_{1z}}
{1+|f_1|^2}\right)_{\bar{z}}+\left(\frac{f_{2z\bar{z}}}{f_{2\bar{z}}}-
\frac{2\bar{f}_2f_{2z}}{1+|f_2|^2}\right)_{\bar{z}}\right\}=0.\nonumber
\end{gather}
Further, there is a natural relationship between the Gauss map (\ref{Gmap}) 
and generalized Weierstrass representation for surfaces, which
defined as follows \cite{KL98a}:
\begin{eqnarray}
X^1+iX^2&=&\int_\Gamma(-\varphi_1\varphi_2dz^{\p}+\psi_1\psi_2d\bar{z}^{\p}),
\nonumber\\
X^1-iX^2&=&\int_\Gamma(\bar{\psi}_1\bar{\psi}_2dz^{\p}-\bar{\varphi}_1
\bar{\varphi}_2d\bar{z}^{\p}),\nonumber\\
X^3+iX^4&=&\int_\Gamma(\varphi_1\bar{\psi}_2dz^{\p}+\psi_1\bar{\varphi}_2
d\bar{z}^{\p}),\nonumber\\
X^3-iX^4&=&\int_\Gamma(\bar{\psi}_1\varphi_2dz^{\p}+\bar{\varphi}_1\psi_2
d\bar{z}^{\p}),\label{W1}
\end{eqnarray}
where
\begin{equation}\ar\label{W2}
\begin{array}{ccc}
\psi_{\alpha z}&=&p\varphi_\alpha,\\
\varphi_{\alpha\bar{z}}&=&-p\psi_{\alpha},
\end{array}\quad\alpha=1,2.
\end{equation}
$\Gamma$ is a contour in complex plane $\C$, $\psi_\alpha,\varphi_\alpha$ 
are complex--valued functions. The formulae (\ref{W1}), (\ref{W2}) define a
conformal immersion of the surface $S^{2,0}$ into the space $\R^{4,0}$.
At this point an induced metric of $S^{2,0}$ has a form \cite{KL98a}:
\[
ds^2=u_1u_2dzd\bar{z},
\]
where $u_\alpha=|\psi_\alpha|^2+|\varphi_\alpha|^2\;(\alpha=1,2)$. Respectively,
gaussian and mean curvature are
\[
K=-\frac{2}{u_1u_2}[\log(u_1u_2)]_{z\bar{z}},\quad H^2=4\frac{|p|^2}
{u_1u_2}.
\]
The generalized Weierstrass representation (\ref{W1}), (\ref{W2}) is related
with the Gauss map (\ref{Gmap}) by means of the following substitutions:
\[
f_1=i\frac{\bar{\psi}_1}{\varphi_1},\quad f_2=-i\frac{\bar{\psi}_2}
{\varphi_2}.
\]

Further, in accordance with theorem \ref{tDir} the Dirac operator on the
surface $S^{2,0}\hookrightarrow M^{4,0}$ has a form
\begin{eqnarray}
D(\psi^+)=(\alpha-iH)\psi^-,\nonumber\\
D(\psi^-)=(\alpha+iH)\psi^+,\nonumber
\end{eqnarray}
where $\alpha=\lambda_1+\lambda_2,\,\beta=2$, and the restricted spinor
field $\Psi_{|S^{2,0}}$ on the surface $S^{2,0}$
according to theorem \ref{timmer}
and relations (\ref{comp}) has the form $\Psi_{|S^{2,0}}=(\psi^+,\psi^-)=
(\phi_1+\phi_3,\phi_4+\phi_2)$, where $\phi_i\in\Om_{0,2}$. Therefore,
in virtue of an inverse Gauss map $g^{-1}$ and the formulae (\ref{comp}) a
Dirac operator on the Riemann surface $S_0$ is equivalent to the following
two systems:
\begin{equation}\ar\label{W3}
\begin{array}{ccc}
\phi_{1z}&=&(\alpha-iH)\phi_4,\\
\phi_{4z^\ast}&=&(\alpha+iH)\phi_1,
\end{array}\quad\begin{array}{ccc}
\phi_{3z}&=&(\alpha-iH)\phi_2,\\
\phi_{2z^\ast}&=&(\alpha+iH)\phi_3,
\end{array}
\end{equation}
where the spinors $\phi_i\in\Om_{0,2}$ on $S_0$ are complex--valued
functions on variables $z,z^\ast$. Let
\begin{eqnarray}
X^1+iX^2&=&\int_\Gamma(-\phi_4\phi_2dz+\phi_1\phi_3dz^\ast),\nonumber\\
X^1-iX^2&=&\int_\Gamma(\phi^\ast_1\phi^\ast_3dz-\phi^\ast_1\phi^ast_2dz^\ast),
\nonumber\\
X^3+iX^4&=&\int_\Gamma(\phi_4\phi^\ast_3dz+\phi_1\phi^\ast_2dz^\ast),\nonumber\\
X^3-iX^4&=&\int_\Gamma(\phi^\ast_1\phi_2dz+\phi^\ast_4\phi_3dz^\ast).
\label{W4}
\end{eqnarray}
Then formulae (\ref{W3}), (\ref{W4}) define a conformal immersion of the
surface $S^{2,0}$ into the space $\R^{4,0}$. At this point an induced
metric has a form
\[
ds^2=(|\phi_1|^2+|\phi_4|^2)(|\phi_3|^2+|\phi_2|^2)dzdz^\ast.
\]
In the case of the parallel spinor field $(\lambda_1=\lambda_2=0)$ 
the formulae (\ref{W3}),
(\ref{W4}) reduce to the generalized Weierstrass representation (\ref{W1}),
(\ref{W2}) if suppose $p=iH$ ¨ $\psi_1=\phi_1,\varphi_1=\phi_4,
\psi_2=\phi_3,\varphi_2=\phi_2$.

Analogously, when the surface $S^{0,2}$ immersed into the Kleinian
manifold $M^{2,2}$ a Dirac operator on $S^{0,2}\hookrightarrow
M^{2,2}$ is defined as follows (theorem \ref{tDir}):
\begin{eqnarray}
D(\psi^+)&=&(-\alpha+iH)\psi^-,\nonumber\\
D(\psi^-)&=&(-\alpha-iH)\psi^+,\nonumber
\end{eqnarray}
where $\alpha=\lambda_1+\lambda_2,\,\beta=-2$, and the restricted spinor
field $\Psi_{|S^{0,2}}=\psi=(\psi^+,\psi^-)$ ($\psi$ is an element of the
minimal left ideal $I_{2,2}\simeq\Om_{2,0}\frac{1}{2}
(1-i\Upsilon_{12})$) is expressed via the spinor $\phi\in\Om_{2,0}$, which
in the matrix representation has a form
\[\ar
\phi=\begin{pmatrix}
\phi_1+\phi_2 & \phi^\ast_4-\phi^\ast_3\\
\phi_3+\phi_4 & \phi^\ast_1-\phi^\ast_2
\end{pmatrix},
\] 
where
\[
\phi_1=a^0+ia^{12},\quad
\phi_2=a^{13}-ia^{23},\quad
\phi_3=a^3-ia^{123},\quad
\phi_4=a^1+ia^2.
\]
At this point $\Psi_{|S^{0,2}}=(\epsilon^+_{02}\psi,\epsilon^-_{02}\psi)=
(\phi_1+\phi_2,\phi_3+\phi_4)$ (theorem \ref{timmer}). Locally, for the
every fiber $\pi^{-1}(x)=T_xM^{2,2}\simeq\R^{2,2}$ there exists a conformal
immersion of the surface $S^{0,2}$ into the space $\R^{2,2}$ defined by the
following formulae:
\begin{eqnarray}
X^1+iX^2&=&\int_\Gamma(\phi_3\phi_4dz+\phi_1\phi_2dz^\ast),\nonumber\\
X^1-iX^2&=&\int_\Gamma(\phi^\ast_1\phi^\ast_2dz+\phi^\ast_3\phi^\ast_4dz^\ast),
\nonumber\\
X^3+iX^4&=&i\int_\Gamma(\phi^\ast_1\phi_4dz+\phi^\ast_3\phi_2dz^\ast),
\nonumber\\
X^3-iX^4&=&-i\int_\Gamma(\phi_3\phi^\ast_2dz+\phi_1\phi^\ast_2dz^\ast),
\nonumber
\end{eqnarray}
where
\[\ar
\begin{array}{ccc}
\phi_{1z}&=&(-\alpha+iH)\phi_3,\\
\phi_{3z^\ast}&=&(-\alpha-iH)\phi_4,
\end{array}\quad
\begin{array}{ccc}
\phi_{2z}&=&(-\alpha+iH)\phi_4,\\
\phi_{4z^\ast}&=&(-\alpha-iH)\phi_2.
\end{array}
\]

Let us consider now a most interesting case (from the viewpoint of physics)
of the immersion of the space--like surface $S^{0,2}$ into the Lorentzian
manifold $M^{1,3}$. According to theorem \ref{tDir} for the Dirac operator on
the surface $S^{0,2}\hookrightarrow M^{1,3}$ we have
\begin{eqnarray}
D(\psi^+)&=&(-\alpha+H)\psi^-,\nonumber\\
D(\psi^-)&=&(-\alpha-H)\psi^+,\nonumber
\end{eqnarray}
where the restricted spinor field $\Psi_{|S^{0,2}}=\psi=(\psi^+,\psi^-)$ is
expressed via the Dirac--Hestenes spinor field $\phi\in\C_2$ by the formulae
(\ref{restrict2}), $\psi$ is the element of the minimal left ideal
$I_{1,3}\simeq\C_2\frac{1}{2}(1+
i\sigma_{12})$ (theorem \ref{timmer}). Further, since the every fiber of
the tangent bundle of $M^{1,3}$ is isomorphic to the Minkowski spacetime,
$\pi(x)=T_xM^{1,3}\simeq\R^{1,3}$,
then there exists (for the every fiber) a conformal immersion of the
surface $S^{0,2}$ into the spacetime
$\R^{1,3}$. This immersion may be defined as follows
\begin{eqnarray}
X^1&=&\frac{1}{2}\int_\Gamma\left[(\phi_1\phi_4+\phi_3\phi_2)dz+
(\phi^\ast_1\phi^\ast_4+\phi^\ast_3\phi^\ast_2)dz^\ast\right],\nonumber\\
X^2&=&\frac{1}{2}\int_\Gamma\left[(\phi_1\phi_4-\phi_3\phi_2)dz+
(\phi^\ast_1\phi^\ast_4-\phi^\ast_3\phi^\ast_2)dz^\ast\right],\nonumber\\
X^3&=&\frac{i}{2}\int_\Gamma\left[(\phi_4\phi_3-\phi_1\phi_2)dz+
(\phi^\ast_4\phi^\ast_3-\phi^\ast_1\phi^\ast_2)dz^\ast\right],\nonumber\\
X^4&=&\frac{1}{2}\int_\Gamma\left[(\phi_4\phi_3+\phi_1\phi_2)dz-
(\phi^\ast_1\phi^\ast_2+\phi^\ast_4\phi^\ast_3)dz^\ast\right],\label{W5}
\end{eqnarray}
where
\begin{equation}\ar\label{W6}
\begin{array}{ccc}
\phi^\ast_{1z}&=&(-\alpha+H)\phi_4,\\
\phi_{4z^\ast}&=&(-\alpha-H)\phi^\ast_1,
\end{array}\quad
\begin{array}{ccc}
\phi^\ast_{3z}&=&-(-\alpha+H)\phi_2,\\
-\phi_{2z^\ast}&=&(-\alpha-H)\phi^\ast_3.
\end{array}
\end{equation}
At this point an induced metric on the surface $S^{0,2}\hookrightarrow\R^{1,3}$
has a form
\[
ds^2=|\phi^\ast_1\phi_2+\phi_4\phi^\ast_3|^2dzdz^\ast.
\]

Further, let $\Psi$ be the parallel spinor field on the manifold $M^{1,3}$, 
then for the system (\ref{W6}) we have
\begin{equation}\ar\label{W7}
\begin{array}{ccc}
\phi^\ast_{1z}&=&H\phi_4,\\
\phi_{4z^\ast}&=&-H\phi^\ast_1,
\end{array}\quad
\begin{array}{ccc}
\phi^\ast_{3z}&=&-H\phi_2,\\
\phi_{2z^\ast}&=&H\phi^\ast_3.
\end{array}
\end{equation}
It is easy to see that the every system (\ref{W7}) coincides with a linear
problem of a modified Veselov--Novikov hierarchy. The first equation of the
mVN--hierarchy has a form \cite{Bog87}
\begin{equation}\label{star}
p_t+p_{zzz}+p_{\bar{z}\bar{z}\bar{z}}+3p_z\omega+3p_{\bar{z}}\bar{\omega}+
\frac{3}{2}p\bar{\omega}_{\bar{z}}+\frac{3}{2}p\omega_z=0,
\end{equation}
where $\omega_{\bar{z}}=(p^2)_z$.\\[0.3cm]
\noindent{\bf Example}.
Suppose now that the surface $S^{0,2}\hookrightarrow\R^{1,3}$ is a surface
of revolution. Then the components of the Dirac--Hestenes spinor field are
defined sa follows
\begin{equation}\ar\label{W8}
\begin{array}{ccc}
\phi^\ast_1&=&r_1(x)\exp(\lambda y),\\
\phi_4&=&s_1(x)\exp(\lambda y),
\end{array}\quad
\begin{array}{ccc}
\phi^\ast_3&=&r_2(x)\exp(\lambda y),\\
\phi_2&=&s_2(x)\exp(\lambda y),
\end{array}
\end{equation}
where $r_i(x),\,s_i(x)$ are real--valued functions, $\lambda\in\C$. The
substitution of (\ref{W8}) into the systems (\ref{W7}), where
\[
\frac{\partial}{\partial z}=\frac{1}{2}\left(\frac{\partial}{\partial x}+
i\frac{\partial}{\partial y}\right),\quad\frac{\partial}{\partial z^\ast}=
\frac{1}{2}\left(\frac{\partial}{\partial x}-i\frac{\partial}{\partial y}
\right),
\]
gives
\begin{equation}\ar\label{W9}
\begin{array}{ccc}
r_{1x}+i\lambda r_1&=&2Hs_1,\\
s_{1x}-i\lambda s_1&=&-2Hr_1,
\end{array}\quad
\begin{array}{ccc}
r_{2x}+i\lambda r_2&=&-2Hs_2,\\
s_{2x}-i\lambda s_2&=&2Hr_2.
\end{array}
\end{equation}
The every system (\ref{W9}) is nothing but a well--known Zakharov--Shabat
system \cite{ZS71}. One--soliton solutions of ZS--system obtained via
the linear Bargmann potentials \cite{Bar49} are well studied 
(see \cite{Lam80}).

Let us show that equations (\ref{W8}) are particular form of a canonical
decomposition of the Dirac--Hestenes spinor field \cite{Hest67}:
\begin{equation}\label{star2}
\phi=r(x)e^{i\beta/2},
\end{equation}
where $r(x)=\sqrt{\rho(x)}R(x)$, $\rho(x)$ is a probability density,
$R(x)\in\spin_+(1,3)$ is a Lorentz rotation, $\beta$ is a so--called
Yvon--Takabayasi angle which defines a duality transformation. Since in our
case the Dirac--Hestenes field is defined on the surface and therefore
depends on two variables, then $r(x)=(x_1,0,0,0)$ and $\beta=(0,x_2,0,0)$,
whilst the spinor (\ref{star2}) depends on four variables $x_1,x_2,x_3,x_4$.
It is obvious that for the spinor field defined on the surface, the variables
$x_3$ and $x_4$ play a role of the deformation parameters. For example,
a dependence on an evolution parameter $x_4=t$ is defined by the standard
procedure of the inverse scattering transform \cite{AS81}. At this point
in the case of the surface of revolution $(p=p(x_1),\,p=u/2)$ the equation
(\ref{star}) reduces to a modified Korteweg--de Vries equation
$u_t=u_{xxx}+3/2u^2u_x$, and a dependence of the potential $u$ on the
parameter $t$ has a form $u=\pm\sech(\mu x-\mu^3t)$, where $\mu$ is a constant
of integration. It allows to express a dependence of fundamental solutions
(Jost functions) of ZS--systems (\ref{W9})
and respectively the Dirac--Hestenes spinor
field (in virtue of (\ref{W8})) on the parameter $t$. Thus, we have a
Dirac--Hestenes spinor field $\phi\in\cl^+_{1,3}$ defined on the surface of
revolution (precisely speaking, solitonic surface of revolution with
reflectionless potential), integrable deformations of which are defined by
the mKdV--hierarchy. In connection with this it should be noted that an
idea of revolution about some fixed axis has deep roots in the electron
theory. For example, Uhlenbeck and Goudsmit in their fundamental paper
\cite{UG25} imagine the electron as a revolving top.
\section*{Acknowledgements}
I am deeply grateful to Prof. H. Baum, to Prof. C. B\"{a}r and 
Prof. P. Lounesto for sending me their interesting papers which are
essentially stimulate me to write down this work.


\begin{thebibliography}{Clifford}
\bibitem[AS81]{AS81} M.J. Ablowitz, H. Segur, {\bf Solitons and the Inverse
Scattering Transform} (SIAM, Philadelphia, 1981). 
\bibitem[Abr89]{Abr89} U. Abresch, {\it Spinor representation of CMC surfaces},
Lecture at Luminy (1989).
\bibitem[Amm98]{Amm98} B. Ammann, {\bf Spin-Strukturen und das Spektrum
des Dirac Operators} (Dissertation Freiburg 1998, Shaker-Verlag, Aachen 1998). 
\bibitem[At71]{At71} M.F. Atiyah, {\it Riemann surfaces and spin structures},
Ann. Scient. Ecole Norm. Sup. {\bf 4}, 47-62 (1971).
\bibitem[Bar49]{Bar49} V. Bargmann, {\it On the connection between phase
shifts and scattering potentials}, Rev. Mod. Phys. {\bf 21}, 488-493 (1949).
\bibitem[B\"{a}r91]{Bar91} C. B\"{a}r, {\bf Das Spektrum von Dirac-Operatoren}
(Dissertation, Bonner Math. Schriften {\bf 217}, 1991). 
\bibitem[B\"{a}r98]{Bar97c} C. B\"{a}r, {\it Extrinsic bounds for
eigenvalues of the Dirac operator}, 
Ann. Glob. Anal. Geom. {\bf 16}, 573--596 (1998).
\bibitem[Bau81]{Bau81} H. Baum, {\bf Spin-Strukturen und Dirac-Operatoren
\"{u}ber pseudoriemannschen Mannigfaltigkeiten} (Teubner, Leipzig, 1981). 
\bibitem[Bau89]{Bau89a} H. Baum, {\it Odd-dimensional Riemannian
manifolds with imaginary Killing spinors}, Ann. Global Anal. Geom. {\bf 7},
141-154 (1989).
\bibitem[BFGK]{BFGK} H. Baum, Th. Friedrich, R. Grunewald, I. Kath,
{\bf Twistors and Killing spinors on Riemannian manifolds} (Teubner-Verlag
Leipzig/Stuttgart 1991).
\bibitem[BD89]{BD89} M. Blau, L. D\c{a}browski, {\it Pin structures on manifolds
quotiented by discrete groups}, J. Geometry and Physics {\bf 6}, 143-157,
(1989).
\bibitem[Bob94]{Bob94} A.I. Bobenko, {\it Surfaces in terms of $2\times 2$
matrices. Old and new integrable cases}, in {\bf Harmonic maps and 
integrable systems} (A. Fordy, J. Wood eds.), 83-127, Vieweg (1994).
\bibitem[Bob99]{Bob99} A.I. Bobenko, {\it Exploring Surfaces through
Methods from the Theory of Integrable Systems. Lectures on the Bonnet
Problem}, preprint SFB 288, N 403, TU-Berlin (1999).
\bibitem[Bog87]{Bog87} L.V. Bogdanov, {\it Veselov--Novikov equation as a
natural two--dimensional generalization of the Korteweg--de Vries equation},
Theor. Math. Phys. {\bf 78}, 309--314 (1987).
\bibitem[BH]{BH} A. Borel, F. Hirzebruch, {\it Characteristic classes and
homogeneous spaces}, Amer. J. Math. {\bf 80}, 458-538, {\bf 81}, 315-382,
{\bf 82}, 491-504 (1958, 1959, 1960).
\bibitem[Che54]{Che54} C. Chevalley, {\bf The Algebraic Theory of Spinors}
(Columbia University Press, New York, 1954).
\bibitem[Che55]{Che55} C. Chevalley, {\it The construction and study of certain
important algebras}, Publications of Mathematical Society of Japan No 1
(Herald Printing, Tokyo, 1955).
\bibitem[Clif78]{cliff} W.K. Clifford, {\it Applications of Grassmann's
extensive algebra}, Amer. J. Math. {\bf 1}, 350, (1878).
\bibitem[Cra85]{Cra85} J. Crawford, {\it On the algebra of Dirac bispinor
densities: Factorization and inversion theorems}, J. Math. Phys. {\bf 26},
1439--1441 (1985).
\bibitem[Cru87]{Cru87} A. Crumeyrolle, {\it The primitive idempotents of the
Clifford algebras and the amorphic spinor fibre bundles}, Reports on Math.
Phys. {\bf 25}, 305-328 (1987).
\bibitem[Cru91]{Cru91} A. Crumeyrolle, {\bf Orthogonal and Symplectic
Clifford Algebras, Spinor Structures} (Kluwer Acad. Publ., Dordrecht, 1991).
\bibitem[Dab88]{Dab88} L. D\c{a}browski, {\bf Group Actions on Spinors}
(Bibliopolis, Naples, 1988).
\bibitem[DLGSC]{DR} C. Doran, A. Lasenby, S. Gull, S. Somaroo, A. Challinor,
{\it Spacetime algebra and electron physics}, Advances in Imaging \&
Electron Physics {\bf 95}, 272-385 (1996).
\bibitem[Eisen]{Eisen} L. P. Eisenhart, {\bf A treatise on the differential
geometry of Curves and Surfaces} (Dover, New York, 1909).
\bibitem[FRO90]{FRO90a} V.L. Figueiredo, W.A. Rodrigues, Jr., E.C. Oliveira,
{\it Covariant, algebraic, and operator spinors}, Int. J. Theor. Phys. {\bf 29},
371-395, (1990).
\bibitem[Fr97]{Fr97} Th. Friedrich, {\bf Dirac-Operatoren in der
Riemannschen Geometrie} (Vieweg-Verlag Braunschweig/Wiesbaden, 1997).
\bibitem[Fr98]{Fr98} Th. Friedrich, {\it On the Spinor Representation
of Surfaces in Euclidean 3-Space},
J. Geometry and Physics {\bf 28}, 143--157 (1998).
\bibitem[GM93]{GM93} P. Ginsparg, G. Moore, {\it Lectures on 2--D gravity
and 2--D string theory}, preprint hep-th/9304011 (1993).
\bibitem[Hae56]{Hae56} A. Haefliger, {\it Sur l'extension du groupe
structural d'un espace fibre}, C.R. Acad. Sci. Paris {\bf 243}, 558-560 (1956).
\bibitem[Hes66]{Hest1} D. Hestenes, {\bf Space--Time Algebra} (Gordon \&
Breach, New York, 1966).
\bibitem[Hes67]{Hest67} D. Hestenes, {\it Real spinor fields}, J. Math.
Phys. {\bf 8}, 798-808, (1967).
\bibitem[Hes76]{Hest76} D. Hestenes, {\it Observables, operators, and
complex numbers in the Dirac theory}, J. Math. Phys. {\bf 16}, 556-571,
(1976).
\bibitem[HO80]{HO80} D.A. Hoffman, R. Osserman, {\it The geometry of
the generalized Gauss map}, Memoirs of the American Mathematical Society 236
(Providence, R.I., 1980).
\bibitem[HO83]{HO83} D.A. Hoffman, R. Osserman, {\it The Gauss map of surfaces
in $R^n$}, J. Differential Geometry {\bf 18}, 733-754 (1983).
\bibitem[HO85]{HO85} D.A. Hoffman, R. Osserman, {\it The Gauss map of surfaces
in $R^3$ and $R^4$}, Proc. London Math. Soc. {\bf 50}, 27-56 (1985).
\bibitem[John80]{John80} D. Johnson, {\it Spin structures and quadratic
forms on surfaces}, J. London Math. Soc. {\bf 22}, 365-373 (1980).
\bibitem[Kar78]{Karo} M. Karoubi, {\bf K-Theory. An Introduction} (Springer-Verlag,
Berlin, 1978).
\bibitem[Kel93]{Kel93} J. Keller, {\it The geometric content of the electron
theory}, Adv. in Appl. Clifford Alg. {\bf 3}(2), 147--200 (1993).
\bibitem[Ken79]{Ken79} K. Kenmotsu, {\it Weierstrass formula for surfaces
of prescribed mean curvature}, Math. Ann. {\bf 245}, 89-99 (1979).
\bibitem[Kon96]{Kon1} B.G. Konopelchenko, {\it Induced surfaces and their
integrable dynamics}, Stud. Appl. Math. {\bf 96}, 9-51, (1996).
\bibitem[Kon98]{Kon2} B.G. Konopelchenko, {\it Weierstrass representation
for surfaces in 4D spaces and their integrable deformations via the DS
hierarchy}, preprint math.DG/9807129, (1998).
\bibitem[KL98a]{KL98a} B.G. Konopelchenko, G. Landolfi, {\it
Generalized Weierstrass representation for surfaces in multidimensional
Riemann spaces}, J. Geometry and Physics {\bf 29}(4), 319--333
(1999).
\bibitem[KL98b]{KonLan2} B.G. Konopelchenko, G. Landolfi, {\it Induced surfaces
and their integrable dynamics. II. Generalized Weierstrass representation
in 4D spaces and deformations via DS hierarchy}, preprint math.DG/9810138,
(1998). 
\bibitem[KT95]{KT95} B.G. Konopelchenko, I.A. Taimanov, {\it Generalized
Weierstrass formulae, soliton equations and Willmore surfaces}, preprint N 187,
Univ. Bochum (1995).
\bibitem[KT96]{KT96} B.G. Konopelchenko, I.A. Taimanov, {\it Constant
mean curvature surfaces via an integrable dynamical system}, J. Phys. A:
Math. Gen. {\bf 29}, 1261-1265 (1996).
\bibitem[KS95]{KS95} R. Kusner, N. Schmitt, {\it The spinor representations
of minimal surfaces}, preprint dg-ga/9512003 (1995). 
\bibitem[KS96]{KS96} R. Kusner, N. Schmitt, {\it The spinor representation
of surfaces in space}, preprint dg-ga/9610005 (1996).
\bibitem[Lam80]{Lam80} G.L. Lamb, Jr., {\bf Elements of soliton theory}
(John Wiley \& Sons, New York, 1980). 
\bibitem[Lou81]{Lou81} P. Lounesto, {\it Scalar Products of Spinors and an
Extension of Brauer-Wall Groups}, Found. Phys. {\bf 11}, 721-740 (1981).
\bibitem[Lou93]{Lou93} P. Lounesto, {\it Clifford algebras and Hestenes
spinors}, Found. Phys. {\bf 23}, 1203-1237 (1993). 
\bibitem[LM89]{LM89} H.B. Lowson, M.-L. Michelsohn, {\bf Spin Geometry}
(Princeton University Press, Princeton 1989).
\bibitem[Mil63]{Mil63} J. Milnor, {\it Spin structures on manifolds},
Enseign. Math. {\bf 9}, 198-203 (1963).
\bibitem[Mil65]{Mil65} J. Milnor, {\it Remarks concerning spin manifolds},
in {\bf Differential and Combinatorial Topology}(Princeton 1965, 55-62
S. Cairns (ed.)). 
\bibitem[Nur96]{Nur96} P. Nurowski, {\it Optical geometries and related
structures}, J. Geometry and Physics {\bf 18}, 335--348 (1996).
\bibitem[Par92]{Par92} J.M. Parra, {\it The Dirac--Hestenes equation and
the algebraic structure of the Minkowski space--time}, XIX ICGTMP,
Salamanca (Spain) July--1992. 
\bibitem[Pen83]{Pen83} R. Penrose, {\it Physical space--time and 
non--realizable CR--structures}, Bull. Amer. Math. Soc. (NS) {\bf 8},
427--448 (1983). 
\bibitem[Port69]{Port} I.R. Porteous, {\bf Topological Geometry} 
(van Nostrand, London, 1969).
\bibitem[Ras55]{Rash} P.K. Rashevskii, {\it The Theory of Spinors} 
(in Russian) Uspekhi Mat. Nauk {\bf 10}, 3--110 (1955); Translation in 
Amer. Math. Soc. Transl. (Ser.2) {\bf 6}, 1 (1957).
\bibitem[Rob61]{Rob61} I. Robinson, {\it Null electromagnetic fields},
J. Math. Phys. {\bf 2}, 290 (1961).
\bibitem[RT86]{RT86} I. Robinson, A. Trautman, {\it Cauchy--Riemann structures
in optical geometry}, Proc. of the Fourth Marcel Grossmann Meeting on General
Relativity, pp. 317--324, 1986.
\bibitem[RSVL]{RSVL} W.A. Rodrigues, Jr., Q.A.G. de Souza, J. Vaz, Jr.,
P. Lounesto, {\it Dirac-Hestenes spinor fields in Riemann-Cartan spacetime},
Int. J. Theor. Phys., {\bf 35}, 1849-1900, (1996).
\bibitem[Roz55]{Roz55} B.A. Rozenfel'd, {\bf Non--Euclidean Geometries}
(Moscow, 1955) [in Russian].
\bibitem[Sul89]{Sul89} D. Sullivan, {\it The spinor representation of 
minimal surfaces}, Notes (1989).
\bibitem[Sym85]{Sym85} A. Sym, {\it Soliton surfaces and their applications},
in: Soliton geometry from spectral problems, Lecture Notes in Physics 239,
Springer, Berlin, 154-231 (1985). 
\bibitem[Tai97a]{Tai97a} I.A. Taimanov, {\it Modified Novikov-Veselov equation
and differential geometry of surfaces}, Trans. Amer. Math. Soc. Ser.2 {\bf 179},
133-159 (1997).
\bibitem[Tai97b]{Tai97b} I.A. Taimanov, {\it Surfaces of revolution in terms
of solitons}, Ann. Glob. Anal. Geom. {\bf 15}, 419-435 (1997).
\bibitem[Tai97c]{Tai97c} I.A. Taimanov, {\it The Weierstrass representation
of closed surfaces in $\R^3$}, preprint SFB 288, N 291, TU-Berlin (1997).
\bibitem[Tai98]{Tai98} I.A. Taimanov, {\it The Weierstrass representation
of spheres in $\R^3$, the Willmore numbers and soliton spheres}, preprint
SFB 288, N 302, TU-Berlin (1998).
\bibitem[UG25]{UG25} G.E. Uhlenbeck, S. Goudsmit, {\it Spinning electrons and
the structure of spectra}, Nature {\bf 117}, 264 (1925).
\bibitem[Var99a]{Var99a} V.V. Varlamov, {\it Generalized Weierstrass
representation for surfaces in terms of Dirac-Hestenes spinor field},
J. Geometry and Physics {\bf 32}(3), 241--251 (1999).
\bibitem[Var99b]{Var99b} V.V. Varlamov, {\it On spinor fields on the
surfaces of revolution}, Proc. Int. Conf. ''Geometrization of Physics IV''
Kazan State University, Kazan, October 4--8, 1999, 248--253.
\bibitem[Var99c]{Var99c} V.V. Varlamov, {\it Fundamental Automorphisms of
Clifford Algebras and an Extension of D\c{a}browski Pin Groups}, 
Hadronic Journal {\bf 22}, 497--533 (1999).
\bibitem[Weier]{Weier} K. Weierstrass, {\it Untersuchungen \"{u}ber die
Fl\"{a}chen, deren mittlere Kr\"{u}mmung \"{u}berall gleich Null ist},
Monatsber. Acad. Wiss. Berlin, 1866, s.612-625.
\bibitem[ZS71]{ZS71} V.E. Zakharov, A.B. Shabat, {\it Exact theory of 
two--focussing and one--dimensional self--modulation of waves in nonlinear
media}, Sov. Phys. -- JETP {\bf 61}, 118--134 (1971).
\end{thebibliography}
\end{document}